\pdfoutput=1

\documentclass[11pt,final,oneside,reqno]{amsart}
\usepackage{mathieu}
\allowdisplaybreaks[3]

\usepackage[color,notref]{showkeys}
\definecolor{refkey}{rgb}{1,0,0}
\definecolor{labelkey}{rgb}{0,0,1}

\hypersetup{pdfpagemode=UseOutlines,colorlinks=false,pdfpagelayout=SinglePage,pdfstartview=FitH,bookmarksopen=true}

\begin{document}
\selectlanguage{english}

\title[Atiyah Classes]{From Atiyah Classes to Homotopy Leibniz Algebras}
\author[ZC]{Zhuo Chen}
\address{Department of Mathematics, Tsinghua University
}
\email{\href{mailto:zchen@math.tsinghua.edu.cn}{zchen@math.tsinghua.edu.cn}}
\author[MS]{Mathieu Sti\'enon}
\address{Department of Mathematics, Penn State University
}
\email{\href{mailto:stienon@math.psu.edu}{stienon@math.psu.edu}}
\author[PX]{Ping Xu}
\address{Department of Mathematics, Penn State University
}
\email{\href{mailto:ping@math.psu.edu}{ping@math.psu.edu}}
\thanks{Research partially supported by NSFC grant 11471179,
the Beijing high education young elite teacher project,
NSA grant H98230-12-1-0234,
and NSF grants  DMS0605725, DMS0801129, DMS1101827.}

\begin{abstract}
A celebrated theorem of Kapranov states that the Atiyah class of the tangent
bundle of a complex manifold $X$ makes $\shiftby{T_X}{-1}$
into a Lie algebra object in $D^+ (X)$, the  bounded below
 derived category of coherent sheaves on $X$.
Furthermore Kapranov proved that, for a K\"ahler manifold $X$,
the Dolbeault resolution $\Omega^{\bullet-1}(T_X^{1, 0})$ of $\shiftby{T_X}{-1}$ is an $L_\infty$ algebra.
In this paper, we prove that Kapranov's theorem holds in much wider generality for vector bundles over Lie pairs.
Given a Lie pair $(L,A)$, i.e.\ a Lie algebroid $L$ together with a Lie subalgebroid $A$,
we define the Atiyah class $\alpha_E$ of an $A$-module $E$ as the obstruction
to the existence of an \emph{$A$-compatible}  $L$-connection on $E$.
We prove that the Atiyah classes $\alpha_{L/A}$ and $\alpha_E$ respectively
make $L/A[-1]$ and $E[-1]$ into a Lie algebra and a Lie algebra module
in the bounded below derived category $D^+(\category{A})$,
where $\category{A}$ is the abelian    category
of left $\enveloppante{A}$-modules and  $\enveloppante{A}$ is the universal
enveloping algebra of $A$.
Moreover, we produce a homotopy Leibniz algebra and a homotopy Leibniz
module stemming from the Atiyah classes of $L/A$ and $E$, and inducing
the aforesaid Lie structures in $D^+(\category{A})$.
\end{abstract}

\maketitle

\tableofcontents

\section*{Introduction}

The Atiyah class of a holomorphic vector bundle $E$ over a
complex manifold $X$, as initially introduced
by Atiyah~\cite{Atiyah}, constitutes the obstruction to the existence of
a holomorphic connection on said holomorphic vector bundle.
It is constructed in the following way. The
vector bundle $\jet{}{E}$ of jets (of order 1) of holomorphic sections of
$E\to X$ fits into the canonical short exact sequence
\[ 0\to T^*_X\otimes E\to\jet{}{E}\to E\to 0 \]
of holomorphic vector bundles
(over the complex manifold $X$).
The Atiyah class of $E\to X$ is the extension class
\[ \alpha_E\in\Ext^1_{X}(E,T^*_X\otimes E)
\cong H^1(X,T^*_X\otimes\End E) \]
of this short exact sequence~\cite{Atiyah,Kapranov}.

In the late 1990's, Rozansky and Witten proposed a construction of
a family of  3-dimensional topological quantum field theories,
indexed by compact (or asymptotically flat) hyper-K\"ahler manifolds \cite{Rozansky_Witten}.
Thus, to each compact hyper-K\"ahler manifold,
the Rozansky-Witten procedure associates
a topological invariant of 3-manifolds.
In subsequent work, Kapranov~\cite{Kapranov} and Kontsevich~\cite{Kontsevich}
revealed the crucial role played by Atiyah classes in the construction of the
Rozansky-Witten  invariants.
In particular, they showed that the hyper-K\"ahler restriction
is unnecessary and that the theory devised by Rozansky and Witten
works for all holomorphic symplectic manifolds~\cite{Kapranov, Kontsevich}.
Kapranov's work highlighted a key fact: the Atiyah class of the tangent
bundle of a complex manifold $X$ yields a morphism
$\shiftby{T_X}{-1}\otimes \shiftby{T_X}{-1} \to \shiftby{T_X}{-1}$
in the bounded below derived category $D^+ (X)$
of coherent sheaves on $X$, which turns $\shiftby{T_X}{-1}$
into a Lie algebra object in $D^+(X)$.
Therefore, Kapranov's approach shone light on many similarities between
the Rozansky-Witten and Chern-Simons theories~\cite{BGRT,BLT}
as stressed by Roberts and Willerton \cite{Willerton}.

Atiyah classes have also enjoyed renewed vigor due to
Kontsevich's seminal work on deformation
quantization \cite{Kontsevich_formality, Kontsevich:motives}.
Kontsevich indicated the   existence of deep ties between the Todd genus
of complex manifolds and the Duflo element in Lie
theory \cite{Kontsevich_formality, Kontsevich:motives, Shoikhet:duflo, MR2816610}.
This discovery inspired several subsequent works on Hochschild (co)homology
and the Hochschild--Kostant--Rosenberg isomorphism for complex manifolds,
by Dolgushev, Tarmarkin \& Tsygan~\cite{MR2565036,arXiv:0807.5117},
 C\u{a}ld\u{a}raru~\cite{Caldararu},
Markarian~\cite{Markarian},
Ramadoss~\cite{Ramadoss},
and Calaque \& Van~den~Bergh~\cite{Calaque} among many others.
In particular, the work of Markarian~\cite{Markarian}
(see also Ramadoss~\cite{Ramadoss})  led to an alternative proof
 of the Hirzebruch--Riemann--Roch theorem  and its variations.

In~\cite{Molino,Molino_topology}, Molino introduced
an Atiyah class for connections ``transverse to a foliation,''
which measures the obstruction to their ``projectability.''
Molino's class has applications in geometry, for instance,
in the study of differential operators on foliated manifolds~\cite{Vaisman}
and in deformation quantization~\cite{MR2223149}.

This paper is the first in a sequence of works
\cite{CSX:CR2014, LSX:CR2012, LSX:2014, VoglarieX:CMP2015, MSX:CR2015} which aim at developing in a general setting
a theory of Atiyah classes and their applications.
Our goal is to explore emerging connections between derived geometry
and classical areas of mathematics such as complex geometry,
foliation theory, Poisson geometry and Lie theory.
The present paper develops  a framework
which encompasses both the original Atiyah class
of holomorphic vector bundles and the Molino class of
real vector bundles foliated over a foliation as special cases.

Lie algebroids  are  the starting point of our approach.
Indeed, holomorphic vector bundles and vector bundles
foliated over a foliation may both be seen as instances
of the concept of module over a Lie algebroid, a straightforward
generalization of the well-known representations of Lie algebras.
Given a Lie algebroid $L$ over a base manifold $M$
with anchor $\rho: L\to TM$, an $L$-connection
on a vector bundle $E\to M$ is a bilinear map
$X\otimes s \mapsto \nabla_{X}s$ from $\Gamma (L)\otimes\Gamma (E)$ to $\Gamma (E)$
satisfying the usual axioms:
 $ \nabla_{fX}s=f\nabla_X s$ and
$\nabla_X (fs)=\rho(X)f\cdot s+f\nabla_X s$ for all $f\in\cinf{M}$.
If the connection is flat, $E$ is said to be a module
 over the Lie algebroid $L$.
When the base is the one-point space,
the $L$-modules are simply Lie algebra modules in the classical sense.
When the base is a complex manifold $X$,
the holomorphic vector bundles over  $X$
are equivalent  to the modules of the complex Lie algebroid $T^{0, 1}_X$
stemming from the complex manifold $X$. Molino's foliated bundles are modules
over the  Lie algebroid structure corresponding to
the foliation of their base.

We introduce a general theory of Atiyah classes of vector bundles over Lie algebroid pairs.
By a Lie algebroid pair $(L,A)$, we mean a Lie algebroid $L$ (over a manifold $M$)
together with a Lie subalgebroid $A$ (over the same base $M$) of $L$.
And by a vector bundle $E$ over the Lie algebroid pair $(L,A)$, we mean a
vector bundle $E$ (over $M$), which is a module over the Lie subalgebroid $A$.
Given such a Lie algebroid pair $(L,A)$ and $A$-module $E$,
we consider the jet bundle $\jet{{\LoverA}}{E}$ (of order 1),
whose smooth sections are the $L$-connections on $E$
extending the (infinitesimal) $A$-action on $E$ in a compatible way.
We prove the following

\begin{introthm}\label{THMA}
The jet bundle $\jet{{\LoverA}}{E}$ is naturally an $A$-module.
It fits in a short exact sequence of $A$-modules
\[ 0\to A^\perp \otimes E\to \jet{{\LoverA}}{E}\to E\to 0 .\]
Here $A^\perp$ denotes the annihilator of $A$ in $L$.
\end{introthm}

We call the extension class
\[ \alpha_E\in \Ext^1_{\category{A}}(E,A^\perp\otimes E)
\cong H^1(A, A^\perp\otimes\End E) \]
of this short exact sequence, the Atiyah class of $E$ because,
when $L=T_X\otimes\CC$ and $A=T_X^{0,1}$ for a complex manifold $X$,
$\jet{{\LoverA}}{E}$ is the space of 1-jets of holomorphic sections of $E\to X$ and $\alpha_E$ is the (classical) Atiyah class of $E\to X$; and, when $L$ is the tangent bundle of a smooth manifold $M$
and $A$ is an integrable distribution on $M$,
$\alpha_E$ is the Molino class of the vector bundle $E$ foliated over $A$.
Geometrically, the Atiyah class can thus be interpreted as the obstruction
to the existence of a compatible $L$-connection on $E$ which extends the
$A$-action with which $E$ is endowed.

It turns out that the Atiyah class introduced in our general context
and the classical Atiyah class of holomorphic vector bundles
enjoy similar algebraic properties.

We denote the abelian category of modules over the universal enveloping algebra
$\enveloppante{A}$ of the Lie algebroid $A$ by the symbol $\category{A}$.
Every vector bundle over $M$ endowed with an $A$-action --- more precisely its space of smooth sections ---
is an object of $\category{A}$.
The bounded below derived category of $\category{A}$ will be denoted by $D^+(\category{A})$.

Given a Lie algebroid pair $(L, A)$, the quotient $L/A$
is naturally an $A$-module. When $L$ is the tangent bundle to a manifold $M$
and $A$ is an integrable distribution on $M$, the $A$-action
on $\bB$ is given by the Bott connection \cite{Bott}.
We consider $L/A$ as a complex concentrated in degree 1 and refer to it as $L/A[-1]$.
We show  that the Atiyah class of $\bB$ makes
$\bB[-1]$ into a Lie algebra object in the derived category $D^+(\category{A})$.
Indeed, being an element of
\[ \begin{split}
\Ext^1_{\category{A}}(\bB, A^\perp\otimes \bB)
\cong\ & \Ext^1_{\category{A}}(\bB\otimes \bB, \bB) \\
\cong\ & \Hom_{D^+(\category{A})}(\bB\otimes\bB,\shiftby{\bB}{1}) \\
\cong\ & \Hom_{D^+(\category{A})}(\shiftby{\bB}{-1}
\otimes\shiftby{\bB}{-1},\shiftby{\bB}{-1})
,\end{split} \]
the Atiyah class $\atiyahclass_{\bB}$ of the $A$-module $\bB$
defines a ``Lie bracket'' on $\shiftby{\bB}{-1}$.
Moreover, if $E$ is an $A$-module, its Atiyah class
\[ \atiyahclass_{E}\in\Ext^1_{\category{A}}(E, A^\perp\otimes E)
\cong\Ext^1_{\category{A}}(\bB\otimes E,E)\cong
\Hom_{D^+(\category{A})}(\shiftby{\bB}{-1}\otimes\shiftby{E}{-1},\shiftby{E}{-1}) \]
defines a ``representation'' on $\shiftby{E}{-1}$ of the ``Lie algebra'' $\shiftby{\bB}{-1}$.
In summary, we prove the following

\begin{introthm}\label{THMB}
Let $(L,A)$ be a Lie algebroid pair with quotient $L/A$. Then
the Atiyah class of $L/A$ turns $\shiftby{\bB}{-1}$ into a Lie algebra object in the derived category $D^+(\category{A})$. Moreover, if $E$ is an $A$-module, then $\shiftby{E}{-1}$
is a module object over the Lie algebra object $\shiftby{\bB}{-1}$
in the derived category $D^+(\category{A})$.
\end{introthm}

The above result suggests that, on the cochain level,
the Atiyah class should define some kind of Lie algebra up to homotopy.
But how does one obtain a cocycle representing the Atiyah class?
Recall that a Dolbeault representative of the Atiyah class
of a holomorphic vector bundle $E\to X$ can be obtained in the following way.
Considering $T\oz_X$ as a complex Lie algebroid, choose a $T\oz_X$-connection
$\nabla\oz$ on $E$. Being a holomorphic vector bundle, $E$ carries a canonical
flat $T\zo_X$-connection $\partialbar$. Adding $\nabla\oz$ and $\partialbar$, we
obtain a $T_X\otimes\CC$-connection $\nabla$ on $E$.
The element $\mathcal{R}\in\OO^{1,1}(\End E)$ defined by
$\mathcal{R}(X^{0,1},Y^{1,0})s=\nabla_{X^{0,1}}\nabla_{Y^{1,0}}s-
\nabla_{Y^{1,0}}\nabla_{X^{0,1}}s-\nabla_{[X^{0,1},Y^{1,0}]}s$
is a Dolbeault 1-cocycle whose cohomology class (which is independent of the
choice of $\nabla\oz$) is the Atiyah class $\alpha_E\in H^{1,1}(X,\End E)
\cong H^1(X,T^*_X\otimes\End E)$.
In the more general setting of vector bundles over Lie algebroid pairs,
the Atiyah class can be recovered from a Lie algebroid connection as follows.
Assume $(L,A)$ is a Lie pair, and $E$ is an $A$-module.
Let $\nabla$ be any $L$-connection on $E$ extending its $A$-action.
The curvature of $\nabla$ is the bundle map $R^\nabla:\wedge^2 L\to\End E$
defined by $R^\nabla(l_1,l_2)=\nabla_{l_1}\nabla_{l_2}-\nabla_{l_2}\nabla_{l_1}
-\nabla_{\lie{l_1}{l_2}}$, for all $l_1, l_2\in\sections{L}$.
Since $E$ is an $A$-module, the restriction of $R^\nabla$
to $\wedge^2 A$ vanishes.
Hence the curvature determines a section
$\atiyahcocycle_E\in\sections{A^*\otimes A^\perp\otimes\End E}$,
which happens to be a $1$-cocycle for the Lie algebroid $A$ with values in the $A$-module
$A^\perp\otimes\End E$.
We prove that the cohomology class $\atiyahclass_E\in H^1(A,A^\perp\otimes\End E)$
of the 1-cocycle $\atiyahcocycle_E$ is the Atiyah class of the $A$-module $E$.

When the $A$-module $E$ is the quotient $\bB$ of the Lie algebroid pair $(L,A)$,
the choice of an $L$-connection $\nabla$ on $\bB$ extending the $A$-action,
yields the Atiyah $1$-cocycle $\atiyahcocycle_{\bB}\in\sections{A^*\otimes(\bB)^*\otimes\End(\bB)}$,
which may be regarded as a bundle map $R_2:\bB\otimes\bB\to A^*\otimes\bB$.
Starting from $R_2$ and a splitting of the short exact sequence of
vector bundles $0\to A\to L\to \bB\to 0$, a sequence $(R_n)_{n=2}^{\infty}$ of bundle maps
\begin{equation*}
R_n:\otimes^n \bB\to\Hom(A,\bB)
\end{equation*}
can  be defined inductively by the relation
$R_{n+1}=\partial^\nabla R_n$, for $n\geq 2$.
Alternatively, $R_n$ can be seen as a section of the vector bundle
$A^*\otimes(\otimes^n (\bB)^*)\otimes\bB$.
Then the graded vector space
$V=\bigoplus_{n=0}^\infty\sections{\wedge^n A^*\otimes\bB}$
can be endowed with a sequence $(\lambda_k)_{k\in\NN}$
of multibrackets $\lambda_k:\otimes^kV\to V$: the unary bracket
$\lambda_1$  is chosen to be the coboundary operator $\partial^A$
of exterior forms on the Lie algebroid $A$ taking values
 in the $A$-module $\bB$,
while all higher order brackets $\lambda_k$ are defined by the relation
\begin{equation*}
\lambda_k(\xi_1\otimes b_1,\cdots,\xi_k\otimes b_k)
=\minuspower{\abs{\xi_1}+\cdots+\abs{\xi_k}}
\xi_1\wedge\cdots\wedge\xi_k\wedge R_k(b_1,\cdots, b_k)
,\end{equation*}
where $b_1,\dots,b_k\in\sections{\bB}$ and $\xi_1,\dots,\xi_k$
are arbitrary homogeneous elements of $\sections{\wedge^\bullet A^*}$.

By an $A$-algebra, we mean a bundle of associative algebras
$\Cp$ over $M$ endowed with an $A$-module structure such that
$\Gamma (A)$ acts by derivations.
For a commutative $A$-algebra $\Cp$, $(\lambda_k)_{k=1}^\infty$
extends in a natural way to the graded vector space $\bigoplus_{n=0}^\infty
\sections{\wedge^n A^*\otimes\bB\otimes \Cp}$.
We prove

\begin{introthm}\label{THMC}
Assume that $(L,A)$ is a Lie pair and $\Cp$ is a commutative $A$-algebra.
When endowed with the sequence of multibrackets $(\lambda_k)_{k\in\NN}$,
the graded vector space
$\sections{\wedge^\bullet A^*\otimes L/A\otimes \Cp}[-1]$
becomes a Leibniz$_\infty$  algebra ---
a natural generalization of Stasheff's $L_\infty$ algebras~\cite{Lada_Stasheff}
first introduced by Loday~\cite{Loday}
in which the requirement that the multibrackets be (skew-)symmetric is dropped.
\newline
If $E$ is an $A$-module, the graded vector space
$\sections{\wedge^\bullet A^*\otimes E \otimes \Cp}[-1]$
becomes a  Leibniz$_\infty$ module over the
Leibniz$_\infty$ algebra
$\sections{\wedge^\bullet A^*\otimes L/A \otimes \Cp}[-1]$.
\newline
As a consequence, $\bigoplus_{i\geq 1}H^{i-1} (A, {\LoverA} \otimes \Cp )$
is a graded Lie algebra and $\bigoplus_{i\geq 1}H^{i-1}(A, E \otimes \Cp)$
 a module over it.
\end{introthm}


We also identify a simple criterion for detecting when this
Leibniz algebra is actually an $L_\infty$ algebra.
This criterion being satisfied when $X$ is a K\"ahler manifold,
$L=T_X\otimes\CC$ and $A=T_X^{0, 1}$, we recover
the $L_\infty$-structure on $\Omega^{0,\bullet-1}(T^{1,0})$
discovered by Kapranov~\cite{Kapranov}.

Recently, we proved that the sequence of multibrackets
$\{\lambda_k\}_{k\in\NN}$ of Theorem~\ref{THMC}
can be tweaked so as to make
$\sections{\wedge^\bullet A^*\otimes L/A}$ an $L_\infty[1]$ algebra
rather than a mere $\Leibniz_\infty[1]$ algebra~\cite{LSX:CR2012,LSX:2014}.
Such an $L_\infty[1]$ algebra should be related to
the $L_\infty$-spaces of Costello~\cite{MR2827826,arXiv:1112.0816}.
The intrinsic meaning of this homotopy algebraic structure
arising from the Atiyah class of a Lie algebroid pairs
is explored in~\cite{LSX:CR2012,LSX:2014,1507.01051}.
Our definition of the Atiyah class can also be
extended to complexes of $A$-modules~\cite{CSX:CR2014}
as was done in~\cite{Markarian}
for complexes of coherent sheaves of $\mathcal{O}_X$-modules.
The universal enveloping algebra of the Lie algebra object $\shiftby{\bB}{-1}$
in the derived category $D^+(\category{A})$ (see Theorem~\ref{THMB})
is described in that same paper~\cite{CSX:CR2014}.
We note that the Atiyah class of Lie algebroid pairs plays a central role
in the construction of new Rozansky--Witten type invariants of 3-manifolds
from symplectic Lie pairs~\cite{VoglarieX:CMP2015}.
In another direction, Atiyah classes were defined for differential graded vector bundles
and these yield homotopy algebraic structures as well~\cite{MSX:CR2015}.

We also would like to point out works of others which are related to the present paper.
Vitagliano studied various homotopy algebra structures associated to foliations
(a special case of Lie algebroid pairs) \cite{MR3277952,MR3300319,MR3313214}.
For the Atiyah class of a DG-module over a dDG-algebra, the reader will want to consult Calaque's work~\cite{Calaque}.
After the first draft of this paper was posted on arXiv, Calaque inferred that, for Lie algebra pairs
$(\mathfrak{d},\mathfrak{g})$, i.e.\ Lie algebroid pairs with the one-point space as base manifold,
the Atiyah class of the quotient $\mathfrak{d}/\mathfrak{g}$ coincides with the class capturing
the obstruction to the ``PBW problem'' studied earlier by Calaque--C\u{a}ld\u{a}raru--Tu~\cite{CCT} (see also~\cite{Grinberg}).
Bordemann gave a nice interpretation of the Calaque--C\u{a}ld\u{a}raru--Tu class
as the obstruction to the existence of invariant connections
on homogeneous spaces~\cite{Bordemann}.
Another recent development is Calaque's beautiful work~\cite{Calaque:PBW}
on the relation between the Atiyah class of the $A$-module $L/A$
with respect to the Lie pair $(L,A)$ and the relative PBW problem
previously solved by C\u{a}ld\u{a}raru--Calaque--Tu~\cite{CCT}.
For more on this topic, we also refer the reader to work of Laurent-Gengoux \& Voglaire~\cite{1507.01051}.
Calaque also pointed out to us that our results should be related to
the obstruction to a relative Hochschild--Kostant--Rosenberg isomorphism for
closed embeddings of algebraic varieties identified by Arinkin \& C\u{a}ld\u{a}raru~\cite{AC}.
This certainly deserves further investigation.
Finally, we would like to mention, in relation to the homotopy algebra
 results of the present paper,
Yu's work on $L_\infty$-algebroids~\cite{Yu}.
Stasheff's work on constrained Poisson algebras \cite{MR940489}
is another interesting result which could
 well be related to the present paper.


\subsection*{Acknowledgments} We would like to express our gratitude to several institutions for their hospitality
while we were working on this project: Penn State University (Chen), Universit\'e Paris~7 (Sti\'enon),
Universit\'e du Luxembourg (Chen and Sti\'enon), Institut des Hautes \'Etudes Scientifiques (Xu),
and Beijing International Center for Mathematical Research (Xu). We would also like to thank  Martin Bordemann,
Paul Bressler, Damien Calaque, Murray Gerstenhaber, Gr\'egory Ginot, Camille Laurent-Gengoux, Boris Shoikhet,
Jim Stasheff, Izu Vaisman and Alan Weinstein for fruitful discussions and useful comments.
Special thanks go to Jim Stasheff for carefully reading the preliminary version of the manuscript.

\section{Preliminaries: connections, modules, Lie pairs, and matched pairs}

Let $M$ be a smooth manifold, let $L\to M$ be a Lie algebroid, and let
$E\xto{\pi}M$ be a vector bundle. The anchor map of $L$ is denoted by $\anchor$.

Recall that the Lie algebroid differential
$d:\Gamma (\wedge^{\bullet}L^*)\to \Gamma ( \wedge^{\bullet+1}L^* )$ is
given by \begin{multline*}
\big(d\mu\big)(x_0,\cdots,x_n)
= \sum_{i=0}^n (-1)^i \rho(x_i) \big( \mu(x_0,\cdots,\widehat{x_i},\cdots,x_n) \big) \\
+ \sum_{i<j} (-1)^{i+j} \mu(\lie{x_i}{x_j},x_0,\cdots,\widehat{x_i},\cdots,\widehat{x_j},\cdots,x_n),
\\\forall \mu\in \Gamma (\wedge^n L^*), x_i\in \Gamma (L), i=0, \cdots, n.
\end{multline*}

The traditional description of a (linear) $L$-connection on $E$ is in terms of a covariant derivative
\[ \sections{L}\times\sections{E}\to\sections{E}: (x,e)\mapsto \nabla_x e \]
characterized by the following two properties:
\begin{gather*}
\nabla_{fx} e=f\nabla_x e , \\
\nabla_x (fe)=\rho(x)f\cdot e+f\cdot\nabla_x e
,\end{gather*}
for all $x\in\sections{L}$, $e\in\sections{E}$, and $f\in\cinf{M}$.

Here, we give three equivalent descriptions of (linear) $L$-connections on $E$:
covariant differential, horizontal lifting, and horizontal distribution.

\begin{definition}
A (linear) $L$-connection on $E$ is a map
$\sections{E}\xto{d^\nabla}\sections{L^*\otimes E}$,
called covariant differential, satisfying the Leibniz rule
\begin{equation*}
d^\nabla (fe)=\rho^* (df) \otimes e + f \cdot d^\nabla e ,
\end{equation*}
 for all $f\in\cinf{M}$ and $e\in\sections{E}$.
\end{definition}

The covariant differential
$\sections{E}\xto{d^\nabla}\sections{L^*\otimes E}$
extends uniquely to a degree 1 operator
\[ \sections{\wedge^{\bullet} L^*\otimes E}\xto{d^\nabla}
\sections{\wedge^{\bullet+1} L^*\otimes E} \]
satisfying the Leibniz rule
\[ d^\nabla (\beta\otimes e)=d\beta \otimes e
+ (-1)^{b}\beta\otimes d^\nabla e ,\]
for all $\beta\in\sections{\wedge^{b}L^*}$ and $e\in\sections{E}$.

\begin{definition}
A (linear) $L$-connection on $E$ is a map $L\times_M E\xto{h} T_E$,
called horizontal lifting, such that the diagram
\[ \xymatrix{ & & L \ar[rr]^{\rho} \ar[rd] & & T_M \ar[ld] \\
L\times_M E \ar[rr]^{h} \ar[rru] \ar[rd] & & T_E \ar[rru]^{\pi_*} \ar[ld] & M & \\
& E \ar[rru]_{\pi} & & & } \]
commutes and its faces
\[ \xymatrix{ L\times_M E \ar[r]^h \ar[d] & T_E \ar[d] \\ E \ar[r]_{\id} & E }
\qquad \text{and} \qquad
\xymatrix{ L\times_M E \ar[r]^h \ar[d] & T_E \ar[d]^{\pi_*} \\ L \ar[r]_{\rho} &
T_M } \]
are vector bundle maps.
\end{definition}

A vector field $X$ on $E$ is said to be projectable onto $M$ if
$\pi(e_1)=\pi(e_2)$ implies $\pi_*(X_{e_1})=\pi_*(X_{e_2})$.
By $\XX_\pi(E)$, we denote the space of vector fields $X$ on $E$ which are
projectable onto $M$ and whose flow $\Phi^X_t:E\to E$ is a vector bundle
automorphism over the flow $\Phi^{\pi_* X}_t:M\to M$ of the projected vector
field $\pi_* X$ on $M$. These vector fields
 are normally called \emph{linear vector fields}
on $E$ (see \cite{MackenzieXu:1997} for details).
 The space $\XX_\pi(E)$ of linear vector fields on $E$ is obviously a module
over the ring $\cinf{M}$.

\begin{definition}
A (linear) $L$-connection on $E$ is a morphism of $\cinf{M}$-modules
$\sections{L}\xto{H}\XX_\pi(E)$, called horizontal distribution,
such that the diagram
\[ \xymatrix{ \sections{L} \ar[rr]^H \ar[rd]_\rho && \XX_\pi(E) \ar[ld]^{\pi_*}
\\ & \XX(M) & } \]
commutes.
\end{definition}

Covariant differential, covariant derivative, horizontal lift, and horizontal
distribution are related to one another by the identities
\begin{gather*}
\nabla_l e=\pairing{d^\nabla e}{l} ; \\
e_*\rho(l_x)-h(l_x,e_x)=\tau_{e_x} (\nabla_l e)_x ; \\
{H(l)}|_{e_x}=h(l_x,e_x) ,
\end{gather*}
$\forall x\in M$, $l\in\sections{L}$, and $e\in\sections{E}$.
Here in the second equation, $\tau_{e_x}$ denotes the canonical isomorphism
between the fiber $E_x$ and its tangent space at the point $e_x$.
This second equation can be rewritten as
\begin{equation}\label{Eqt:handnabla}
h(l_x,e_x)f_{\nu}=\pairing{\nabla_{l_x}\nu}{e_x}
,\end{equation}
where $f_{\nu}$ denotes the fiberwise linear function on $E$ determined by $\nu\in\sections{E^*}$.

The following assertions are equivalent:
\begin{gather*}
\nabla_{l_1}\nabla_{l_2}-\nabla_{l_2}\nabla_{l_1}=\nabla_{\lie{l_1}{l_2}} ;\\
H(\lie{l_1}{l_2})=\lie{H(l_1)}{H(l_2)} .
\end{gather*}
When they are satisfied for all $l_1,l_2\in\sections{L}$, the connection is said to be flat.
An $L$-module is a vector bundle $E \to M$ endowed with a flat (linear) $L$-connection.
A flat (linear) $L$-connection will also be called an $L$-action or $L$-representation.
When the $L$-connection on $E$ is flat, $(d^\nabla)^2=0$ and
$(\sections{\wedge^\bullet L\otimes E}, d^\nabla)$
is a cochain complex, whose cohomology groups $H\graded(L, E)$ are
the so-called Lie algebroid cohomology groups
of $L$ with values in $E$.

By a \emph{Lie algebroid pair}, or simply a \emph{Lie pair} $(L,A)$,
we mean a Lie algebroid $L$ and a Lie subalgebroid $A$ of $L$ over the same
base manifold $M$.

\begin{proposition}\label{Prop:quotientmodule}
The quotient $L/A$ of a Lie pair $(L,A)$ is an $A$-module; the action of $A$ on $L/A$ is defined by
\[ \nabla_a \bigl(\quotientmapLB(l)\bigr)=\quotientmapLB(\lie{a}{l}), \quad
\forall a\in\sections{A}, l\in\sections{L} ,\]
where $\quotientmapLB$ denotes the projection $L\to L/A$.
Being dual to $L/A$, the annihilator $A^\perp$ of $A$ in $L^*$ is also an $A$-module.
\end{proposition}

Assume now that $A$ and $B$ are two Lie subalgebroids of a Lie algebroid $L$
such that $L$ and $A\oplus B$ are isomorphic as vector bundles. Then $L/A\isomorphism B$
is naturally an $A$-module while $L/B\isomorphism A$ is naturally a $B$-module.
The Lie algebroids $A$ and $B$ are said to form a matched pair.

\begin{definition}[\cite{Lu,Mokri,Mackenzie}]
Two (real or complex) Lie algebroids $A$ and $B$ over the same base manifold $M$
and with respective anchors $\rho_A$ and $\rho_B$ are said to form a matched pair
if there exists an action $\nabla$ of $A$ on $B$ and an action $\anadelta$ of $B$ on $A$
such that the identities
\begin{gather*}
\lie{\anchor_A(X)}{\anchor_B(Y)} = -\anchor_A\big(\anadelta_Y X\big)+\anchor_B\big(\nabla_X Y\big) , \\
\nabla_X\lie{Y_1}{Y_2} = \lie{\nabla_X Y_1}{Y_2} +
\lie{Y_1}{\nabla_X Y_2} + \nabla_{\anadelta_{Y_2} X}Y_1 -
\nabla_{\anadelta_{Y_1} X} Y_2 , \\
\anadelta_Y\lie{X_1}{X_2} = \lie{\anadelta_Y X_1}{X_2} +
\lie{X_1}{\anadelta_Y X_2} + \anadelta_{\nabla_{X_2} Y} X_1 -
\anadelta_{\nabla_{X_1}Y}X_2 ,
\end{gather*}
hold for all $X_1,X_2,X\in\sections{A}$ and $Y_1,Y_2,Y\in\sections{B}$.
\end{definition}

\begin{proposition}[\cite{Mokri,Mackenzie}]\label{thm:5.2}
Given a matched pair $(A,B)$ of Lie algebroids, there is a Lie
algebroid structure $A\bowtie B$ on the direct sum vector bundle
$A\oplus B$, with anchor \[ X\oplus Y\mapsto\rho_A(X)+\rho_B(Y) \] and bracket
\begin{equation*}
\lie{X_1\oplus Y_1}{X_2\oplus Y_2}
= \big(\lie{X_1}{X_2} + \anadelta_{Y_1}X_2 - \anadelta_{Y_2}X_1 \big)
\oplus \big( \lie{Y_1}{Y_2} + \nabla_{X_1}Y_2 - \nabla_{X_2}Y_1 \big) .
\end{equation*}
Conversely, if $A\oplus B$ has a Lie algebroid structure for
which $A\oplus 0$ and $0\oplus B$ are Lie subalgebroids, then the
representations $\nabla$ and $\anadelta$ defined by
\[ \lie{X\oplus 0}{0\oplus Y} = -\anadelta_Y X\oplus\nabla_X Y \]
endow the couple $(A,B)$ with a structure of matched pair.
\end{proposition}

\begin{example}
A Lie algebra is a Lie algebroid whose base manifold is the one-point space.
If the direct sum $\mfg\oplus\mfg^*$ of a vector space $\mfg$ and its dual $\mfg^*$
is endowed with a Lie algebra structure such that the direct summands $\mfg$ and $\mfg^*$ are Lie subalgebras and
\[ [X,\alpha]=\coadjoint{X}\alpha - \coadjoint{\alpha}X
,\qquad \forall X\in\mfg, \alpha\in\mfg^* ,\]
the pair $(\mfg,\mfg^*)$ is said to be a Lie bialgebra.
Lie bialgebras are instances of matched pairs of Lie algebroids.
\end{example}

\begin{example}\label{t01t10}
Let $X$ be a complex manifold. Then $(T^{0,1}_X,T^{1,0}_X)$ is a
matched pair, where the actions are given by
\[ \nabla_{X^{0,1}}X^{1,0}=\pr^{1,0}\lie{X^{0,1}}{X^{1,0}}
\qquad \text{ and } \qquad
\anadelta_{X^{1,0}}X^{0,1}=\pr^{0,1}\lie{X^{1,0}}{X^{0,1}}, \] for all
$X^{0,1}\in\sections{T\zo_X}$ and $X^{1,0}\in\sections{T\oz_X}$.
Hence $T^{0,1}_X\bowtie T^{1,0}_X$ and $T_X\otimes\CC$
are isomorphic as complex Lie algebroids.
More generally, given a holomorphic Lie algebroid $A$,
the pair $(A^{0,1},A^{1,0})$
is a matched pair of Lie algebroids and $A^{0,1}\bowtie A^{1,0}$
is isomorphic, as a complex Lie algebroid, to $A\otimes\CC$ \cite{LSX}.
\end{example}

\begin{example}
Let $D$ be an integrable distribution on a smooth manifold $M$.
Then $D$ is a Lie subalgebroid of $T_X$, and
the normal bundle $T_X/D$ is canonically a $D$-module.
The $D$-action on $T_X/D$ is usually called Bott connection \cite{Bott}.
Moreover, if $\mathcal{F}_1$ and $\mathcal{F}_2$ are two transversal foliations
on a smooth manifold $M$, the corresponding tangent distributions
$T_{\mathcal{F}_1}$ and $T_{\mathcal{F}_2}$ constitute a matched pair of Lie
algebroids with $T_{\mathcal{F}_1}\bowtie T_{\mathcal{F}_2}\isomorphism T_X$.
\end{example}

\begin{example}
Let $G$ be a Poisson Lie group and let $(P,\pi)$ be a Poisson $G$-space,
i.e. a Poisson manifold $P$ together
with a $G$-action such that the action map
 $G\times P\to P$ is  a Poisson map.
According to Lu \cite{Lu}, $A=(T^*P)_\pi$ and $B=P\rtimes\mathfrak{g}$ form a matched
pair of Lie algebroids.
\end{example}

\begin{remark}
A matched pair of Lie algebroids $L=A\bowtie B$ can be seen as a Lie pair $(L,A)$
together with a splitting $j:B\to L$ of the short exact sequence $0\to A \to L \to B \to 0$,
whose image $j(B)$ happens to be a Lie subalgebroid of $L$.
\end{remark}

\section{Atiyah classes}

\subsection{Prelude: holomorphic connections}

The Atiyah class of a holomorphic vector bundle $E$ over a complex manifold $X$
is the obstruction class to the existence of
a holomorphic (linear) connection. It is constructed in the following way. The
vector bundle $\jet{}{E}$ of jets (of order 1) of holomorphic sections of
$E\to X$ fits into the canonical short exact sequence of holomorphic vector bundles
\[ 0\to T^*_X\otimes E\to\jet{}{E}\to E\to 0 \]
over the complex manifold $X$.
The Atiyah class of $E\to X$ is the extension class
\[ \alpha_E\in\Ext^1_{X}(E,T^*_X\otimes E) \]
of this short exact sequence~\cite{Atiyah,Kapranov}.

There are canonical isomorphisms between the abelian groups
$\Ext^1_{X}(E,T^*_X\otimes E)$ and $\Hom_{D^+(X)}(T_X\otimes E,E[1])$,
the sheaf cohomology group $H^1(X,T^*_X\otimes\End E)$
and the Dolbeault cohomology group $H^{1, 1}(X,\End E)$.
A Dolbeault representative of the Atiyah class can be obtained in the following way.
Considering $T\oz_X$ as a complex Lie algebroid, choose a $T\oz_X$-connection
$\nabla\oz$ on $E$. Being a holomorphic vector bundle, $E$ carries a canonical
flat $T\zo_X$-connection $\partialbar$. Adding $\nabla\oz$ and $\partialbar$, we
obtain a $T_X\otimes\CC$-connection $\nabla$ on $E$.
The element $\mathcal{R}\in\OO^{1,1}(\End E)$ defined by
\begin{equation}\label{Eqt:RinOmega11}
\mathcal{R}(X^{0,1}, Y^{1,0})s=\nabla_{X^{0,1}}\nabla_{Y^{1,0}}s-
\nabla_{Y^{1,0}}\nabla_{X^{0,1}}s-\nabla_{[X^{0,1},Y^{1,0}]}s \end{equation}
is a Dolbeault 1-cocycle whose cohomology class (which is independent of the
choice of $\nabla\oz$) is the Atiyah class $\alpha_E\in H^{1, 1}(X,\End E)$.

\subsection{Existence of $A$-compatible $L$-connections}

\subsubsection{Extension of an $A$-action to a compatible $L$-connection}

Throughout this section, $(L,A)$ is a Lie pair and $E$ is an $A$-module.
The symbols $\sheaf{E}$ and $\sheaf{L}$ will denote the sheaves on $M$ defined by
\[ \sheaf{E}(U)=\set{e\in\sections{U;E}\text{ s.t. }\nabla_{a}e=0,
\forall a \in \sections{U;A} } ,\] and
\[ \sheaf{L}(U)=\set{l\in\sections{U;L}\text{ s.t.
}\lie{a}{l}\in\sections{U;A},\forall a\in\sections{U;A}} ,\]
where $U$ denotes an arbitrary open subset of $M$.

\begin{lemma}\label{Lem:nabla12difference}
Given an $A$-module $E$, there always exists an $L$-connection on $E$ extending
the given $A$-connection. Moreover, if $\nabla^1$ and $\nabla^2$ are two such
extensions, then $d^{\nabla^2}-d^{\nabla^1}\in\sections{A^\perp\otimes\End E}$,
where $A^\perp$ denotes the annihilator of $A$ in $L^*$.
\end{lemma}

\begin{proof}
Choose a subbundle $B$ of $L$ such that $L=A\oplus B$
and a $T_M$-connection $\nabla^{(T_M)}$ on $E$ --- this is always possible.
Then extend $A$-connection $\nabla^{(A)}$ to an $L$-connection $\nabla^{(L)}$ on $E$ by setting
\[ \nabla^{(L)}_{a+b}= \nabla^{(A)}_a+\nabla^{(T_M)}_{\rho(b)} ,\]
where $\rho$ denotes the anchor map $L\to T_M$.
The difference $l \mapsto \nabla^1_l-\nabla^2_l$ of two such extensions $\nabla^1$ and $\nabla^2$
is a bundle map $L\to\End E$, which vanishes on $A$.
\end{proof}

\begin{definition}
An $L$-connection $\nabla$ on $E$ is said to be
$A$-compatible,  if
\begin{enumerate}
\item
  it extends the given $A$-action on $E$, and
\item
 it satisfies
\[ \nabla_a \nabla_l-\nabla_l\nabla_a=\nabla_{\lie{a}{l}}, \quad
\forall a\in\sections{A}, l\in\sections{L} .\]
\end{enumerate}
\end{definition}

\begin{proposition}
Let $\nabla$ be an $L$-connection on $E$ extending its $A$-action.
Provided that the sheaf of smooth sections of $E$ is isomorphic to
$C^{\infty}_M\otimes_{\RR}\sheaf{E}$ and the sheaf of smooth sections
of $L$ is isomorphic to $C^{\infty}_M\otimes_{\RR}\sheaf{L}$,
the $L$-connection $\nabla$ is $A$-compatible if and only if
$\nabla_{\sheaf{L}}\sheaf{E}\subset\sheaf{E}$.
\end{proposition}

\begin{proof}
For any $a\in\sections{U;A}$, $l\in\sheaf{L}(U)$ and $e\in\sheaf{E}(U)$, we have
$\lie{a}{l}\in\sections{U;A}$, $\nabla_a e=0$, and $\nabla_{\lie{a}{e}}=0$
so that, if $\nabla$ is $A$-compatible, we obtain
\[ \nabla_a \nabla_l e=\nabla_a \nabla_l e-\nabla_l \nabla_a e-\nabla_{\ba{a}{l}} e=0 .\]
Hence $\nabla_l e\in\sheaf{E}(U)$. This proves that
$\nabla_{\sheaf{L}}\sheaf{E}\subset\sheaf{E}$.
Conversely, if $\nabla_{\sheaf{L}}\sheaf{E}\subset\sheaf{E}$, then
\[ \bigl(\nabla_a\nabla_{f\cdot l}-\nabla_{f\cdot l}\nabla_a
-\nabla_{\ba{a}{ {f\cdot l}}}\bigr)(g\cdot e) = fg\cdot \nabla_{a}\nabla_{l}e =0 ,\]
for all $a\in\sections{U;A}$, $f,g\in C^{\infty}_U$, $l\in\sheaf{L}(U)$ and $e\in\sheaf{E}(U)$.
Since $\sections{U;E}=C^{\infty}_U\otimes_{\RR}\sheaf{E}$
and  $\sections{U;L}=C^{\infty}_M\otimes_{\RR}\sheaf{L}$,
it follows that $\nabla$ is $A$-compatible.
\end{proof}

\begin{remark}
\label{rmk:14}
Given a matched pair of Lie algebroids $(A,B)$ and an $A$-module $E$,
consider the Lie algebroid $L=A\bowtie B$.
An $A$-compatible $L$-connection on $E$
determines a $B$-connection on $E$ satisfying
\[ \nabla_a\nabla_b e-\nabla_b\nabla_a e=\nabla_{\lie{a}{b}} e,
\quad \forall a\in\sections{A}, b\in\sections{B}, e\in\sections{E} .\] The converse is also true.
\end{remark}

\subsubsection{Atiyah class: obstruction to compatibility}

Assume that $(L,A)$ is a Lie pair,
$E$ is an $A$-module, and $\nabla$ is an $L$-connection on $E$ extending its $A$-action.
The curvature of $\nabla$ is the bundle map
$R^\nabla:\wedge^2 L\to\End E $ defined by
\begin{equation*}
R^\nabla(l_1,l_2)=\nabla_{l_1}\nabla_{l_2}-\nabla_{l_2}\nabla_{l_1}
-\nabla_{\lie{l_1}{l_2}}, \quad \forall~ l_1, l_2\in\sections{L}.
\end{equation*}
Since $E$ is an $A$-module, its restriction to $\wedge^2 A$ vanishes.
Hence the curvature induces a section
$\atiyahcocycle_E\in\sections{A^*\otimes A^\perp\otimes\End E}$
or, equivalently, a bundle map $\atiyahcocycle_E:A\otimes(L/A)\to\End E$
given by
\begin{equation}\label{Eqt:defnatyahE}
\atiyahcocycle_E\big(a;q(l)\big)
=R^\nabla(a,l)=\nabla_{a}\nabla_{l}-\nabla_{l}\nabla_{a}-\nabla_{\lie{a}{l}},\quad
\forall a\in\sections{A},l\in\sections{L}
.\end{equation}
The $L$-connection $\nabla$ is compatible with the $A$-module structure of $E$
if and only if $\atiyahcocycle_E=0$.

\begin{theorem}\label{Thm:atiyahclass}
\begin{enumerate}
\item The section $\atiyahcocycle_E$ of $A^*\otimes A^\perp\otimes\End E $
is a $1$-cocycle for the Lie algebroid $A$ with values in the $A$-module
$A^\perp\otimes\End E $.
\item The cohomology class $\atiyahclass_E\in H^1(A, A^\perp\otimes\End E)$
of the $1$-cocycle $\atiyahcocycle_E$
does not depend on the choice of $L$-connections extending the $A$-action.
\item The Atiyah class $\atiyahclass_E$ of $E$ vanishes if and only if
there exists an $A$-compatible $L$-connection on $E$.
\end{enumerate}
\end{theorem}

\begin{proof}
We use the symbol $\partial^A$ to denote the covariant differential associated to
the action of the Lie algebroid $A$ on $A^\perp\otimes\End E$.
\newline
(a) The second Bianchi identity states that
$d^\nabla R^\nabla:\wedge^3 L\to\End E$ is identically zero.
Thus, for any $a_1,a_2\in\sections{A}$, $l\in\sections{L}$, we have
\begin{align*}
0 =& (d^\nabla R^\nabla)(a_1,a_2,l) \\
=& \nabla_{a_1}(R^\nabla(a_2,l))-\nabla_{a_2}(R^\nabla(a_1,l))+\nabla_{l}(R^\nabla(a_1,a_2)) \\
& -R^\nabla(\ba{a_1}{a_2},l)+R^\nabla(\ba{a_1}{l},a_2)-R^\nabla(\ba{a_2}{l},a_1) \\
=& \nabla_{a_1}\big(\atiyahcocycle_E(a_2;\quotientmapLB(l))\big)
-\nabla_{a_2}\big(\atiyahcocycle_E(a_1;\quotientmapLB(l))\big) \\
& -\atiyahcocycle_E(\ba{a_1}{a_2};\quotientmapLB(l))
-\atiyahcocycle_E(a_2;\nabla_{a_1}\quotientmapLB(l))
+\atiyahcocycle_E(a_1;\nabla_{a_2}\quotientmapLB(l)) \\
=& \Big(\nabla_{a_1}\big(\atiyahcocycle_E(a_2;\quotientmapLB(l))\big)
-\atiyahcocycle_E\big(a_2;\nabla_{a_1}\quotientmapLB(l)\big)\Big) \\
& -\Big(\nabla_{a_2}\big(\atiyahcocycle_E(a_1;\quotientmapLB(l))\big)
-\atiyahcocycle_E\big(a_1;\nabla_{a_2}\quotientmapLB(l)\big)\Big)
-\atiyahcocycle_E\big(\ba{a_1}{a_2};\quotientmapLB(l)\big) \\
=& \big(\partial^A\atiyahcocycle_E\big)(a_1,a_2;\quotientmapLB(l))
.\end{align*}
Therefore $\partial^A\atiyahcocycle_E=0$.
\newline
(b) By Lemma~\ref{Lem:nabla12difference},
if $\nabla^1$ and $\nabla^2$ are two $L$-connections that extend the $A$-action,
then $\nabla^1_l-\nabla^2_l=\phi( l )$ for some $\phi\in\sections{A^\perp\otimes\End E}$,  and
\begin{align*}
& R^{\nabla_1}_E(a;\quotientmapLB(l))\cdot e-R^{\nabla_2}_E(a;\quotientmapLB(l))\cdot e \\
&\qquad = \nabla_a(\nabla^1_l-\nabla^2_l)e-(\nabla^1_l-\nabla^2_l)\nabla_a e
-(\nabla^1_{\ba{a}{l}}-\nabla^2_{\ba{a}{l}})e \\
&\qquad = \nabla_a\big(\phi(l)\cdot e\big)-\phi(l)\cdot(\nabla_a e)-\phi(\ba{a}{l})e \\
&\qquad = \big(\partial^A\phi\big)(a;l)\cdot e
.\end{align*}
So $R^{\nabla_1}_E-R^{\nabla_2}_E=\partial^A\phi$.
\newline
(c) It is clear that $\atiyahcocycle_E$ vanishes if and only if $\nabla$ is $A$-compatible.
Now, if $\atiyahcocycle_E=\partial^A\phi$ for some $\phi\in\sections{A^\perp\otimes\End E}$,
set $\nabla'=\nabla-\phi$. Then $R^{\nabla'}_E=0$, which implies that $\nabla'$ is
$A$-compatible.
\end{proof}

We call $\atiyahcocycle_E$ the  \emph{Atiyah cocycle}
associated with the $L$-connection
$\nabla$ that extends the $A$-module structure of $E$,
and  its  corresponding cohomology class
$\atiyahclass_E\in H^1(A, A^\perp\otimes\End E)$
the  \emph{Atiyah class} of the $A$-module $E$.

\begin{remark}\label{rmk:dDG}
When the Lie pair $(L, A)$ is a  matched pair of Lie algebroids,
i.e.\ $L=A\bowtie B$, our definition of Atiyah class is a special case of
the Atiyah class of a dDG algebra developed
by Calaque and Van den Bergh~\cite{Calaque}.
Hence in the matched pair case, Theorem~\ref{Thm:atiyahclass} (a)-(b)
is a consequence of Lemma~8.2.4 in~\cite{Calaque}.
\end{remark}

\begin{example}
Let $X$ be a complex manifold, and $E$ a holomorphic vector bundle over $X$.
Then $A=T_X^{0,1}$ and $B=T_X^{1, 0}$ form a matched pair of Lie algebroids and
$L=A\bowtie B$ is isomorphic to $T_X\otimes\CC$.
Moreover $E$ is an $A$-module \cite{LSX}. It is simple to see that holomorphic
$T_X$-connections on $E$ are equivalent to $L$-connections on $E$ compatible
with the $A$-action (as well as to $A$-compatible $B$-connections on $E$
--- see Remark~\ref{rmk:14}).
In this case, the Atiyah cocycle is exactly the Dolbeault 1-cocycle $\mathcal{R}$
defined by Equation~\eqref{Eqt:RinOmega11}.
\end{example}

\begin{example}
A holomorphic Lie algebroid $K$ over a complex manifold $X$ yields
a matched pair of complex Lie algebroids  $(T_X\zo,K\oz)$ \cite{LSX}.
The Atiyah class of the $T_X\zo$-module $K\oz$ is the Atiyah class for $K$
studied extensively  by Calaque and Van den Bergh in~\cite{Calaque}.
\end{example}

\begin{example}
In~\cite{Molino}, Molino introduced an Atiyah class for connections
``transversal to a foliation,'' which measures the obstruction to their ``projectability.''
Although not phrased in the language of Lie algebroids,
his construction is a special case of ours.
Here $L$ is the tangent bundle $T_M$,
$A$ is the tangent bundle to a foliation $\mathcal{F}$ of $M$,
and the $A$-module $E$ is a vector bundle on $M$ foliated over $\mathcal{F}$.
A transversal connection is an $L$-connection on $E$ which extends the $A$-action.
It is said to be projectable precisely if it is $A$-compatible, i.e.\ if it is preserved
by parallel transport along any path tangent to $\mathcal{F}$.
\end{example}

\begin{example}\label{Ex:dhliealgebras}
Let $\mathfrak{g}$ be a Lie subalgebra of a Lie algebra $\mathfrak{d}$.
Given a $\mathfrak{g}$-module $E$ (i.e.\ a Lie algebra morphism
$\boldsymbol{A}:\mathfrak{g}\to\End E$), and a $\mathfrak{d}$-connection on $E$ extending it
(i.e.\ a linear map $\boldsymbol{L}:\mathfrak{d}\to\End E$ whose restriction to $\mathfrak{g}$ is $\boldsymbol{A}$),
the  Atiyah class is the element in the Chevalley-Eilenberg cohomology group
$H^1(\mathfrak{g},\mathfrak{g}^\perp\otimes \End(E))$ determined by
$\partial^{\mathfrak{g}} \boldsymbol{L}$.
(The symbol $\partial^{\mathfrak{g}}$ denotes the Chevalley-Eilenberg coboundary of
$\mathfrak{d}^*\otimes\End(E)$-valued $\mathfrak{g}$-cochains.)
Here $\boldsymbol{L}$ is considered as an element in $\mathfrak{d}^*\otimes\End(E)$,
which is, in general, not in $\mathfrak{g}^\perp\otimes \End(E)$.
Hence, in general, $\partial^{\mathfrak{g}} \boldsymbol{L}$ does not vanish in
$H^1(\mathfrak{g},\mathfrak{g}^\perp\otimes \End(E))$.
\end{example}

The following example is due to  Calaque-C\u{a}ld\u{a}raru-Tu \cite{CCT}.

\begin{example}
Consider the Lie algebra $\mathfrak{sl}_2(\CC)$ and its standard basis
\[ h= \begin{pmatrix} 1 & 0 \\ 0 & -1 \end{pmatrix} ,\qquad
e= \begin{pmatrix} 0 & 1 \\ 0 & 0 \end{pmatrix} ,\qquad
f= \begin{pmatrix} 0 & 0 \\ 1 & 0 \end{pmatrix} .\]
We have \[ [e,f]=h, \qquad [h,e]=2e,  \qquad [h,f]=-2f .\]
Together, the matrices $h$ and $e$ generate the Lie subalgebra $\mathfrak{g}$
of $2\times 2$ traceless upper triangular matrices.
We identify the quotient $\mathfrak{sl}_2(\CC)/\mathfrak{g}$
to the nilpotent Lie subalgebra $\mathfrak{n}$ generated by $f$.
Note that $\mathfrak{g}$ and $\mathfrak{n}$ form a matched pair of Lie algebras with sum
$\mathfrak{g}\oplus\mathfrak{n}=\mathfrak{sl}_2(\CC)$.
The bilinear map $\theta:\mathfrak{n}\otimes\mathfrak{n}\to\mathfrak{n}$ defined by
$\theta(f,f)=f$ is a generator of the one-dimensional $\mathfrak{g}$-module
$\mathfrak{g}^\perp\otimes\End(\mathfrak{n})\isomorphism
\Hom(\mathfrak{n}\otimes\mathfrak{n},\mathfrak{n})$.
The action of $\mathfrak{g}$ on $\Hom(\mathfrak{n}\otimes\mathfrak{n},\mathfrak{n})$
is given by $h\cdot\theta=2\theta$ and $e\cdot \theta=0$.
One checks that the degree 1 cohomology
$H^1(\mathfrak{g}, \mathfrak{g}^\perp\otimes\End(\mathfrak{n}))$
is a one-dimensional vector space generated by the Atiyah class
$\alpha_{\mathfrak{n}}$ of the $\mathfrak{g}$-module $\mathfrak{n}$.
\end{example}

\subsection{Functoriality}

Let $M$ and $N$ be smooth manifolds, $f:N\to M$  a smooth map,
$A$  a Lie algebroid over $M$ with anchor $\rho:A\to TM$,
and $E$  a smooth vector bundle over $M$.

Let $f^* E$ denote the pullback of $E$ through $f$, i.e.\ the
 fibered product of $N$ and $E$ over $M$:
\[ \xymatrix{ f^*E \ar[d] \ar[r] & E \ar[d] \\ N \ar[r]_{f} & M. } \]

If the anchor $\rho$ and the differential of $f$ are transversal
(i.e.\ $f_*(TN)+\rho(A)=TM|_N$), we can consider the fibered product $f^\star A$
of $TN$ and $A$ over $TM$:
\[ \xymatrix{ f^\star A \ar[d] \ar[r] & A \ar[d]^{\rho} \\ TN \ar[r]_{f_*} & TM. } \]
Note that $\rho$ and $f_*$ are  automatically
transversal when $f$ is a surjective submersion or
when the Lie algebroid $A$ is transitive.
It is clear that  $f^\star A $ is a vector bundle over $N$. However, note that
 $f^\star A \neq f^* A$.
The fiber of $f^\star A$ over a point $n\in N$ is
\[ (f^\star A)_n=\left\{(x,a)\in T_nN\oplus A_{f(n)}\middle| f_*(x)=\rho(a) \right\} .\]
The Lie algebroid structure on $A$ induces a Lie algebroid structure on $f^\star A\to N$;
its anchor is the projection $f^\star A\to TN$ and its bracket is given by
\[ \lie{(x_1,a_1)}{(x_2,a_2)}= (\lie{x_1}{x_2},\lie{a_1}{a_2}) ,\]
for any $x_1,x_2\in\XX(N)$ and $a_1,a_2\in\sections{A}$ such that $f_*(x_1)=\rho(a_1)$
and  $f_*(x_2)=\rho(a_2)$ (see \cite{HMackenzie} for details).

\begin{proposition}
\label{pro:pullback}
Let $A$ be a Lie algebroid over $M$ and
 $f:N\to M$ a smooth map whose differential $f_*:TN\to TM$
is transversal to the anchor of $A$. Then
\begin{enumerate}
\item If $E$ is a module over $A$, then $f^*E$ is a module over $f^\star A$.
\item The map $f$ induces a natural  homomorphism
\[ f^\dagger:H^\bullet \big(A, E\big)\to
H^\bullet \big(f^\star A, f^*E\big) .\]
\end{enumerate}
\end{proposition}
\begin{proof}
The first assertion is easily proved if one thinks of Lie algebroid modules in terms of horizontal lifting.
The second assertion follows from a direct verification.
\end{proof}

The following proposition is immediate.

\begin{proposition}
\label{pro:pullbackLiepair}
If $(L,A)$ is a Lie pair over $M$, and $f:N\to M$  a smooth map
whose differential $f_*:TN\to TM$ is transversal to
 the anchor of $A$,
 then $(f^\star L,f^\star A)$ is a Lie pair over $N$.
\end{proposition}

Given a Lie pair $(L,A)$ over a smooth manifold $M$
and a smooth map $f:N\to M$ whose differential $f_*:TN\to TM$
is transversal to the anchor of $A$
(otherwise $f^\star A$ and $f^\star L$ could be singular),
there is a canonical morphism of vector bundles $f^\star L \to f^* (L/A)$ over $N$:
\[ f^\star L \ni (x_n,a_{f(n)})\mapsto a_{f(n)}+A_{f(n)} \in f^*\left(\frac{L}{A}\right) ,\]
whose kernel is exactly $f^\star A$.
In other words, we have an exact sequence of vector bundles
\[ 0 \to f^\star A \to f^\star L \to f^*(L/A) .\]
Therefore $f^\star L / f^\star A$ can be seen as a vector subbundle of $f^*(L/A)$.

\begin{lemma}
Under the hypothesis of Proposition \ref{pro:pullbackLiepair},
the inclusion
 \[ I: \frac{f^\star L}{f^\star A}\to f^*\left(\frac{L}{A}\right) \]
 intertwines the $ f^\star A$-module structures of $(f^\star L)/(f^\star A)$
and $f^* (L/A)$.
\end{lemma}

Dualizing it, as  a consequence, we obtain the epimorphism of vector bundles
\[ I^\dagger: f^*(A^\perp)\to(f^\star A)^\perp ,\]
which is a morphism of $ f^\star A$-modules.
Note that, when $f$ is a surjective submersion, $I$ is surjective and thus
both $I$ and $I^\dagger$ are isomorphisms of $ f^\star A$-modules.

We are now ready to state the main result in this subsection.

\begin{theorem}
\label{thm:pullback}
Let $(L,A)$ be  a Lie pair over $M$, and $f:N\to M$  a smooth map
whose differential $f_*:TN\to TM$ is transversal to the anchor of $A$.
Assume that $E$ is an $A$-module. Then
the composition of homomorphisms
\[ H^1\big(A,A^\perp\otimes\End E\big)\xto{f^\dagger}
H^1\big(f^\star A,f^*(A^\perp\otimes\End E)\big)\xto{I^\dagger}
H^1\big(f^\star A,(f^\star A)^\perp\otimes\End (f^* E)\big) \]
maps the Atiyah class of $E$ relative to the Lie pair $(L,A)$
onto the Atiyah class of $f^* E$ relative to the Lie pair $(f^\star L,f^\star A)$:
\[ (I^\dagger\rond f^\dagger) (\alpha_E)= \alpha_{f^* E}.\]
\end{theorem}

\subsection{Scalar Atiyah classes and Todd class}

Let  $(L,A)$ be  a Lie pair. We define the scalar Atiyah classes \cite{Atiyah}
of an $A$-module $E$ by
\[ c_k (E):=\frac{1}{k!} \left(\frac{i}{2\pi}\right)^k \trace\big(\alpha_E^k\big)
\in H^k(A,\wedge^k A^\perp) .\]
Here $\alpha_E^k$ denotes the image of
$\alpha_E\otimes\cdots\otimes\alpha_E$
under the natural map
\[ H^1(A,A^\perp\otimes\End E) \times \cdots \times H^1(A,A^\perp\otimes\End E)
\to H^k(A,\wedge^k A^\perp \otimes\End E) \]
induced by the composition in $\End E$ and the wedge product in $\wedge^\bullet A^\perp$.

\begin{remark}
If $E$ is a holomorphic vector bundle over a compact K\"ahler manifold $X$,
the natural inclusion of $H^k(X,\Omega^k)$ into $H^{2k}(X,\CC)$
maps the scalar Atiyah classes of $E$ relative to the Lie pair
$(L=T_X\otimes\CC,A=T\zo_X)$ to the Chern classes of $E$.
\end{remark}

The Todd class of the $A$-module $E$ relative to the Lie pair $(L,A)$ is the cohomology class
\[ \Todd(E) =\det\left(\frac{\alpha_E}{1-e^{-\alpha_E}}\right)
\in H^{\bullet}(A, \wedge^\bullet A^\perp ) .\]

The following propositions can be verified directly.

\begin{proposition}
Let  $(L,A)$ be  a Lie pair and let $E_1$, $E_2$ be $A$-modules. Then
\[ \Todd(E_1\oplus E_2)=\Todd E_1 \cdot \Todd E_2 .\]
\end{proposition}

\begin{proposition}
Under the hypothesis of Theorem~\ref{thm:pullback}, we have
\begin{gather*}
c_k(f^*E) = (I^\dagger\rond f^\dagger) (c_k E)\in H^k\big(f^\star A, \wedge^k (f^\star A)^\perp \big) \\
\Todd(f^*E)=(I^\dagger\rond f^\dagger )(\Todd E)\in H^\bullet\big(f^\star A, \wedge^\bullet (f^\star A)^\perp \big)
\end{gather*}
\end{proposition}

\subsection{Jet short exact sequence}

\subsubsection{The jet bundle $\jet{\LoverA}{E}$}

Let $M$ be a smooth manifold, let $L\to M$ be a Lie algebroid, and let
$E\xto{\pi}M$ be a vector bundle.

An $L$-jet (of order 1) on $E$ (at $e_x\in E$) is a linear map
$L_{\pi(e_x)}\xto{\phi}T_{e_x}E$ such that the diagram
\[ \xymatrix{ L_{\pi({e_x})} \ar[rr]^\phi \ar[rd]_\rho &&
T_{e_x} E \ar[ld]^{\pi_{*{e_x}}} \\ & T_{\pi({e_x})}M & } \]
commutes.
The jet space $\jet{L}{E}$ is the manifold whose points are $L$-jets on $E$.
It is a vector bundle over $X$: the projection $\jet{L}{E}\to M$ maps
$(L_{\pi({e_x})}\xto{\phi} T_{e_x} E)$ to $\pi({e_x})$.
It fits into the short exact sequence of vector bundles over $M$:
\begin{equation}\label{exact0}
\xymatrix{ 0 \ar[r] & L^*\otimes E \ar[r]^{\hat{f}} & \jet{L}{E}\ar[r]^{\hat{g}} & E
\ar[r] & 0 .}
\end{equation}
The surjection $\hat{g}$ maps $(L_{\pi({e_x})}\xto{\phi} T_{e_x} E)$
 to ${e_x}$,
while the injection $\hat{f}$ maps $(L_x\xto{\psi} E_x)$ to
$(L_x\xto{\rho\oplus\psi}T_x M\oplus E_x\isomorphism T_{0_x}E)$.

The following  result is straightforward.

\begin{proposition}
A splitting $s:E\to\jet{L}{E}$ of the short exact sequence of vector bundles
\eqref{exact0} determines a (linear) $L$-connection on $E$. The
converse is also true.
\end{proposition}

In general, there is no canonical choice of splitting for \eqref{exact0}.
However, the induced short exact sequence
\begin{equation}
\label{eq:short-section}
\xymatrix{ 0 \ar[r] & \sections{L^*\otimes E} \ar[r]^{\hat{f}_\sharp} &
\sections{\jet{L}{E}} \ar[r]^{\hat{g}_\sharp} & \sections{E} \ar[r] & 0 }
\end{equation}
at the level of spaces of smooth sections splits canonically:
if $e$ is a section of $E$, then $\goulot{e}:=e_*\rond\rho$ is a
section of $\jet{L}{E}$
such that $\hat{g}_\sharp(\goulot{e})=e$.

We note that the covariant differential
$d^\nabla:\sections{E}\to\sections{L^*\otimes E}$
associated to a splitting $s:E\to\jet{L}{E}$ of the short exact
sequence \eqref{exact0} is given by
\begin{equation*}
\hat{f}_\sharp(d^\nabla e)=\goulot{e}-s_\sharp(e) ,\quad \forall e\in\sections{E}
.\end{equation*}
		
Now assume $A$ is a Lie subalgebroid of $L$ and $E$ is an $A$-module.
The symbol $h$ will denote the horizontal lifting associated to the $A$-action on $E$.

An $h$-extending $L$-jet (of order 1) on $E$ is a linear map
$L_{\pi({e_x})}\xto{\phi} T_{e_x} E$ such that the diagram
\[ \xymatrix{ & A_{\pi({e_x})} \ar[dl] \ar[dr]^{h(-,{e_x})} & \\
L_{\pi({e_x})} \ar[rr]^\phi \ar[rd]_\rho & & T_{e_x} E \ar[ld]^{\pi_{*{e_x}}} \\
& T_{\pi({e_x})}M & } \]
commutes.
The jet space $\jet{\LoverA}{E}$ is the manifold whose points are $h$-extending
$L$-jets on $E$.
It is a vector bundle over $M$: the projection $\jet{{\LoverA}}{E}\to M$ maps
$(L_{\pi({e_x})}\xto{\phi} T_{e_x} E)$ to $\pi({e_x})$.
\begin{example}
When $E$ is a holomorphic vector bundle over a complex manifold  $X$, $A=T_X^{0,1}$ and
$L=T_X\otimes\CC$, the jet bundle $\jet{{\LoverA}}{E}$
is simply the bundle of jets (of order 1) of holomorphic sections of $E$.
\end{example}
Consider the surjective morphism of vector bundles $\breve{g}:\jet{{\LoverA}}{E}\to E$, which maps
$(L_{\pi({e_x})}\xto{\phi} T_{e_x} E)$ to ${e_x}$. Since $T_{0_x}E$ is canonically isomorphic
to $T_x M\oplus E_x$, the kernel of $\breve{g}$ can be identified naturally
with the subbundle
$K$ of $L^*\otimes E \to M$ consisting of all linear maps
 $(L_x\xto{\psi}E_x)$ which satisfy
\[ h(a_x,0_x)=\rho(a_x)+\psi(a_x), \quad\forall x\in M, a_x\in A_x .\]
Since the $A$-connection $h$ on $E$ is linear, $h(a_x,0_x)$ must be the image of $\rho(a_x)$ under the differential of the
zero section $M\xto{0}E$. Therefore, a linear map $(L_x\xto{\psi}E_x)$
 is an element of $K$ if and only if $\psi(a_x)=0$
for all $a_x\in A$, so that $K\cong A^\perp\otimes E$.
Hence we obtain the short exact sequence of vector bundles
\begin{equation}\label{exact2}
\xymatrix{ 0 \ar[r] & A^\perp\otimes E \ar[r]^{\breve{f}}
& \jet{{\LoverA}}{E} \ar[r]^{\breve{g}} & E \ar[r] & 0 ,}
\end{equation}
where the injection $\breve{f}$ maps $(L_x\xto{\psi} E_x)$ to the jet
\[ L_x\to T_{0_x}E\isomorphism T_x M\oplus E_x \qquad  l_x\mapsto\rho(l_x)+\psi(l_x) .\]

In general, there is no canonical choice of splitting for \eqref{exact2}.

\begin{proposition}
\label{trashcan}
A splitting $s:E\to\jet{{\LoverA}}{E}$ of the short exact sequence of vector bundles
\eqref{exact2} determines a (linear) $L$-connection on $E$ extending the $A$-action $h$.
The converse is also true.
\end{proposition}

Obviously, we have the commutative diagram with exact rows
\[ \xymatrix{
0 \ar[r] & A^\perp\otimes E \ar[r]^{\breve{f}} \ar[d] & \jet{L/A}{E} \ar[r]^{\breve{g}} \ar[d] &
E \ar[r] \ar[d]^{\id} & 0 \\
0  \ar[r] & L^*\otimes E \ar[r]^{\hat{f}} & \jet{L}{E} \ar[r]^{\hat{g}} & E \ar[r] & 0 ,} \]
all of whose vertical arrows denote inclusions.

\subsubsection{An equivalent description of the jet bundle}

\begin{proposition}\label{Eqt:equivalentlanguageofjet}
An $L$-jet (of order 1) on $E$ extending the $A$-action $\nabla$
is a pair $(D_x,e_x)$ consisting of a linear map
$D_x:\sections{E^*}\to L_x^*$ and a point $e_x$ in the fiber of $E$ over $x\in M$,
satisfying
\begin{gather}
\pairing{D_x(\varepsilon)}{a_x}=\pairing{\nabla_{a_x}\varepsilon}{e_x}
\quad(\text{or equivalently}\quad D_x(\varepsilon)=\pairing{d^\nabla\varepsilon}{e_x})
; \label{Eqt:jetdefn1} \\
D_x(f\varepsilon)=f(x)\cdot D_x(\varepsilon)+\pairing{\varepsilon_x}{e_x}\cdot \rho^*(df)
, \label{Eqt:jetdefn2}
\end{gather}
for all $a_x\in A_x$, $\varepsilon\in\sections{E^*}$, and $f\in\cinf{M}$.
\end{proposition}

\begin{proof}
Given such a pair $(D_x,e_x)$, each $l_x\in L_x$ determines uniquely a tangent vector
$\tau_x\in T_{e_x}E$ through the relations
\[ \tau_x(\pi^* f)=\rho(l_x)f=\pairing{\rho^*(df)}{l_x} \qquad \text{and} \qquad
\tau_x (f_{\varepsilon})=\pairing{D_x(\varepsilon)}{l_x} ,\]
where $\varepsilon\in\sections{E^*}$, $f\in\cinf{M}$, and $f_{\varepsilon}\in\cinf{E}$
is the fiberwise linear function associated to $\varepsilon$:
\[f_{\varepsilon}(e_x)=\pairing{\varepsilon_x}{e_x} .\]
Let $\phi_D:L_x\to T_{e_x}E$ be the map $l_x\mapsto\tau_x$.
Clearly, $\phi_D$  is linear and satisfies $\pi_*\circ\phi_D=\rho$.
Moreover, $\phi_D$ is an extension of the $A$-action since
\begin{equation*}
\big(\phi_D(a_x)\big)(f_\varepsilon)=\pairing{D_x(\varepsilon)}{a_x}
=\pairing{\nabla_{a_x}\varepsilon}{e_x}=h(a_x,e_x)(f_\varepsilon).
\end{equation*}
Here we have made use of Equation~\eqref{Eqt:handnabla}.
Hence $\phi_D\in\big(\jet{{\LoverA}}E\big)_x$ and $\breve{g}(\phi_D)=e_x$.
Conversely, given an element $\phi:L_x\to T_{e_x}E$ of $\big(\jet{{\LoverA}}E\big)_x$
that projects to $e_x$ under $\breve{g}$, we can define a linear map
$D^\phi_x:\sections{E^*}\to L^*_x$ by the relation
\[ \pairing{D^\phi_x(\varepsilon)}{l_x}=\big(\phi(l_x)\big)(f_{\varepsilon}) .\]
It is straightforward to check that $(D^\phi_x,e_x)$ satisfies \eqref{Eqt:jetdefn1} and \eqref{Eqt:jetdefn2}.
\end{proof}

\begin{remark}
The surjection $\breve{g}:\jet{{\LoverA}}{E}\to E$ in \eqref{exact2}
maps the 1-jet $(D_x,e_x)$ to $e_x$.
The injection $\breve{f}$ in \eqref{exact2} maps $\psi\in (A^\perp\otimes E)_x$
to $(\psi^\dagger_x,0_x)\in\big(\jet{{\LoverA}}{E}\big)_x$, where
$0_x$ is the zero vector of $E_x$ and
$\psi^\dagger_x:\sections{E^*}\to L_x^*$ is the linear map defined by
\begin{equation}
 \pairing{\psi^\dagger_x(\varepsilon)}{l_x}=\pairing{\varepsilon_x}{\psi(l_x)}
,\quad\forall l_x\in L_x,\varepsilon\in\sections{E^*} .
\label{eq:paris}
\end{equation}
Here $\psi$ is considered as a linear map $L_x\to E_x$ whose kernel contains $A_x$.
\end{remark}

\subsubsection{The jet bundle as an $A$-module}

The jet bundle $\jet{{\LoverA}}E$ can be naturally endowed with an $A$-action.

In the language of Proposition~\ref{Eqt:equivalentlanguageofjet}, a section of $\jet{{\LoverA}}E\to M$
consists of a section $e$ of $E\to M$ and an $\RR$-linear map $D:\sections{E^*}\to\sections{L^*}$ satisfying
\begin{gather*}
\pairing{D(\varepsilon)}{a}=\pairing{\nabla_a\varepsilon}{e} ;\\
D(f\varepsilon)=f\cdot D(\varepsilon)+\pairing{\varepsilon}{e}\cdot \rho^*(df)
,\end{gather*}
for all $f\in\cinf{M}$, $\varepsilon\in\sections{E^*}$, and $a\in\sections{A}$.

\begin{proposition}\label{Prop:sequenceofAmodules}
\begin{enumerate}
\item The jet bundle $\jet{{\LoverA}}E$ is a module over $A$; the covariant derivative
\[ \sections{A}\times\sections{\jet{{\LoverA}}E}\to\sections{\jet{{\LoverA}}E} \]
is given by
\[ \nabla_a (D,e)=(\nabla_a D,\nabla_a e) ,\]
where the $\RR$-linear map $\nabla_a D:\sections{E^*}\to\sections{L^*}$
is defined as follows:
\begin{equation}\label{badweather}
\pairing{\big(\nabla_a D\big)(\varepsilon)}{l}=\rho(a)\pairing{D(\varepsilon)}{l}
-\pairing{D(\nabla_a\varepsilon)}{l}-\pairing{D(\varepsilon)}{\lie{a}{l}}
.\end{equation}
\item Diagram~\eqref{exact2} is a short exact sequence of $A$-modules.
\end{enumerate}
\end{proposition}

\begin{proof}
(a) Checking that \eqref{badweather} determines a connection is straightforward.
The flatness of the $A$-connection on $\jet{{\LoverA}}E$ is a consequence of the flatness of the $A$-connection on $E$.
\newline (b) By definition, $\breve{g}$ is a morphism of $A$-modules.
Let us check that $\breve{f}$ is also a morphism of $A$-modules.
For any $\psi\in\sections{A^\perp\otimes E}$, we have
\begin{multline*}
\pairing{(\nabla_a \psi^\dagger)(\varepsilon)}{l}
= \anchor(a)\pairing{\psi^\dagger(\varepsilon)}{l}
-\pairing{\psi^\dagger(\nabla_a\varepsilon)}{l}
-\pairing{\psi^\dagger(\varepsilon)}{\lie{a}{l}} \\
= \anchor(a)\pairing{\psi(l)}{\varepsilon}
-\pairing{\psi(l)}{\nabla_a\varepsilon}
-\pairing{\psi(\lie{a}{l})}{\varepsilon}
= \pairing{(\nabla_a\psi)(l)}{\varepsilon}
= \pairing{(\nabla_a\psi)^\dagger(\varepsilon)}{l}
.\qedhere\end{multline*}
Here the map $\psi^\dagger$ is defined by Equation \eqref{eq:paris}.
\end{proof}

\subsubsection{Alternative description of the $A$-action on $\jet{{\LoverA}}E$}

In this section, $B$  denotes the quotient $L/A$ of the Lie pair $(L,A)$.

The proof of the following lemma is a tedious computation, which we omit.

\begin{lemma}\label{nonlinearity}
The splitting $\goulot{}:\sections{E}\to\sections{\jet{L}{E}}$
of the short exact sequence \eqref{eq:short-section}
 is not $\cinf{M}$-linear.
For every $e\in\sections{E}$ and $f\in\cinf{M}$, we have
\[ \goulot{(f\cdot e)}-f\cdot\goulot{e}=\hat{f}_\sharp\big(\rho^*(df)\otimes e\big) .\]
\end{lemma}

In general, $\goulot{e}$ need not be a section of $\jet{L/A}{E}$.
Nevertheless, fixing a splitting of the short exact sequence of vector bundles
\[ \xymatrix{0 \ar[r] & A \ar[r]^i & L \ar[r]^q & B \ar[r] & 0 } ,\]
i.e.\ a pair of maps $j:B\to L$ and $p:L\to A$ such that
$q\rond j=\id_B$, $p\rond i=\id_A$ and $i\rond p+j\rond q=\id_L$:
\[ \xymatrix{0 \ar@<.5ex>[r] & A \ar@<.5ex>[r]^i \ar@<.5ex>[l] & L \ar@<.5ex>[r]^q \ar@<.5ex>[l]^p & B
\ar@<.5ex>[r] \ar@<.5ex>[l]^j & 0 \ar@<.5ex>[l] } ,\]
naturally determines a splitting $\bravo{}:\sections{E}\to\sections{\jet{L/A}{E}}$
of the short exact sequence of spaces of smooth sections
\begin{equation*}
\xymatrix{ 0 \ar[r] & \sections{A^\perp\otimes E} \ar[r]^{\breve{f}_\sharp} &
\sections{\jet{L/A}{E}} \ar[r]^{\breve{g}_\sharp} & \sections{E} \ar[r] & 0 }
\end{equation*}
induced by \eqref{exact2}.
The image of $x\in M$ under the section $\bravo{e}$ of $\jet{L/A}{E}$
associated to a section $e$ of $E$  by the splitting $\bravo{}$ is the 1-jet
\[ L_x\ni l_x \xmapsto{(\bravo{e})_x}
h(p(l_x),e_x)+e_{*x}\big(\rho\rond j\rond q(l_x)\big) \in T_{e_x}E .\]
It is not difficult to see that $\breve{g}_\sharp(\bravo{e})=e$.

The proof of the following lemma is a tedious computation, which we omit.
\begin{lemma}\label{nonlinearityagain}
The splitting $\bravo{}:\sections{E}\to\sections{\jet{L/A}{E}}$ is not $\cinf{M}$-linear.
For every $e\in\sections{E}$ and $f\in\cinf{M}$, we have
\[ \bravo{(f\cdot e)}-f\cdot\bravo{e}=\breve{f}_\sharp
\Big(\big(l\mapsto\rho(j\rond q(l))f\big)\otimes e\Big) .\]
\end{lemma}

Since both $E$ and $A^\perp$ are modules over $A$, so is $A^\perp\otimes E$:
\[ \pairing{\nabla_a (\lambda\otimes e)}{l}=
\rho(a)\big(\lambda(l)\big)\cdot e-\lambda(\lie{a}{l})\cdot e
+\lambda(l)\cdot\nabla_a e ,\]
where $a\in\sections{A}$, $\lambda\in\sections{A^\perp}$, $e\in\sections{E}$, and $l\in\sections{L}$.

\begin{remark}
Since $\sections{\jet{L/A}{E}}$ is the direct sum of
$\breve{f}_\sharp ( \sections{A^\perp\otimes E})$ and
$\bravo{}(\sections{E})$,
to define an $A$-module structure on the jet bundle $\jet{L/A}{E}$,
it suffices to define the $A$-action on
$\breve{f}_\sharp ( \sections{A^\perp\otimes E})$ and
$\bravo{}(\sections{E})$. Naively, one might expect to define
such an $A$-action by setting
\begin{gather}
\nabla_a\big(\breve{f}_\sharp(\lambda\otimes e)\big)
=\breve{f}_\sharp\big(\nabla_a(\lambda\otimes e)\big) \label{eq:london} ;\\
\nabla_a(\bravo{e})=\bravo{(\nabla_a e)} \label{eq:rome}
.\end{gather}
However, according to a  tedious computation, we have

\[ \bravo{(\nabla_{f\cdot a}e)-\nabla_{f\cdot a}(\bravo{e})}
= f\cdot\big(\bravo{(\nabla_a e)}-\nabla_a(\bravo{e})\big)
+ \breve{f}_\sharp\Big(\big(l\mapsto\rho(j\rond q(l))f\big)\otimes\nabla_a e\Big) \]
and
\[ \bravo{\big(\nabla_a (f\cdot e)\big)}-\nabla_{a}\big(\bravo{(f\cdot e)}\big)
= f\cdot\big(\bravo{(\nabla_a e)}-\nabla_a(\bravo{e})\big)
+ \breve{f}_\sharp\Big(\big(l\mapsto\rho(i\rond p\lie{j\rond q(l)}{a})f\big)\otimes e\Big) ,\]
for any $a\in\sections{A}$, $e\in\sections{E}$, and $f\in\cinf{M}$.
This means that Equation \eqref{eq:rome} must be modified.
\end{remark}

To modify Equation \eqref{eq:rome}, let us introduce some new notations.
Given $a\in\sections{A}$ and $e\in\sections{E}$,
define $\Theta(a,e)\in\sections{A^\perp\otimes E}$ by the relation
\[ \pairing{\Theta(a,e)}{l}=\nabla_{i\rond p\lie{j\rond q(l)}{a}}e,
\qquad \forall l\in\sections{L} .\]

The proof of the following Lemma is a tedious computation, which we omit.

\begin{lemma}\label{coucou}
 For any $f\in\cinf{M}$,
\begin{gather*}
\Theta(f\cdot a,e)-f\cdot\Theta(a,e)
=\big(l\mapsto\rho(j\rond q(l))f\big)\otimes\nabla_a e ;\\
\Theta(a,f\cdot e)-f\cdot\Theta(a,e)
=\big(l\mapsto\rho(i\rond p\lie{j\rond q(l)}{a})f\big)\otimes e
.\end{gather*}
\end{lemma}

This suggests the following proposition.

\begin{proposition}\label{22times}
\begin{enumerate}
\item The jet bundle $\jet{L/A}{E}\to M$ is an $A$-module:
the flat $A$-connection on $\jet{L/A}{E}$ is given by
\begin{gather}
\nabla_a \big(\breve{f}_\sharp(\lambda\otimes e)\big)
=\breve{f}_\sharp\big(\nabla_a (\lambda\otimes e)\big) \label{gather:tp1} \\
\nabla_a (\bravo{e})=\bravo{(\nabla_a e)}-\breve{f}_\sharp\big(\Theta(a,e)\big) , \label{gather:tp2}
\end{gather}
for all $a\in\sections{A}$, $\lambda\in\sections{A^\perp}$,
and $e\in\sections{E}$.
\item Diagram~\eqref{exact2} is a short exact sequence of $A$-modules.
\end{enumerate}
\end{proposition}

\begin{proof}[Sketch of proof]
(a) It follows from Lemmas~\ref{nonlinearityagain} and~\ref{coucou} that what we have defined
is indeed a covariant derivative.
Flatness follows from the Jacobi identity of  the Lie
algebroid $L$, and the flatness of the $A$-connections
on $E$ and $A^\perp\otimes E$.
\newline
(b) It suffices to check that $\breve{g}$ is a morphism of $A$-modules.
We have \[ \breve{g}_\sharp\big(\nabla_a (\bravo{e})\big)=
\breve{g}_\sharp\big(\bravo{(\nabla_a e)}\big) - \breve{g}_\sharp \breve{f}_\sharp\big(\Theta(a,e)\big)
= \nabla_a e = \nabla_a \big( \breve{g}_\sharp (\bravo{e})\big) .\qedhere \]
\end{proof}


\begin{remark}
Observe that, for a matched pair $L=A\bowtie B$,
 $\Theta(a,e)\in\sections{B^*\otimes E}$
is given by the simple formula:
\[ \pairing{\Theta(a,e)}{b}=\nabla_{\anadelta_b a}e .\]
\end{remark}

\begin{proposition}
The $A$-actions on $\jet{\LoverA}E$ defined in Propositions~\ref{Prop:sequenceofAmodules} and~\ref{22times} are identical.
\end{proposition}

\begin{proof}
The $\RR$-linear map $D^{\bravo{e}}:\sections{E^*}\to\sections{L^*}$
determined by the section $\bravo{e}$ of $\jet{\LoverA}E\to M$
(as explained in the proof of Proposition~\ref{Eqt:equivalentlanguageofjet}) satisfies
\begin{gather*}
\pairing{D^{\bravo{e}}(\varepsilon)}{a}=\pairing{\nabla_a\varepsilon}{e} ;\\
\pairing{D^{\bravo{e}}(\varepsilon)}{j(b)}=\Big(\bravo{e}\big(j(b)\big)\Big)(f_{\varepsilon})
=\Big(e_*\rond\rho\big(j(b)\big)\Big)(f_{\varepsilon})=\rho\big(j(b)\big)\pairing{\varepsilon}{e}
,\end{gather*}
for all $e\in\sections{E}$, $\varepsilon\in\sections{E^*}$, $a\in\sections{A}$, and $b\in\sections{B}$.
The image of $D^{\bravo{e}}$ under the action of $a\in\sections{A}$ is the $\RR$-linear map
$\nabla_a D^{\bravo{e}}:\sections{E^*}\to\sections{L^*}$
defined by \eqref{badweather}.
We have
\begin{align*}
& \pairing{\big(\nabla_{a'} D^{\bravo{e}}\big)(\varepsilon)}{a} \\
&\qquad = \rho(a')\pairing{D^{\bravo{e}}(\varepsilon)}{a}-\pairing{D^{\bravo{e}}(\nabla_{a'}\varepsilon)}{a}
-\pairing{D^{\bravo{e}}(\varepsilon)}{\lie{a'}{a}} \\
&\qquad = \rho(a')\pairing{\nabla_a \varepsilon}{e}-\pairing{\nabla_a\nabla_{a'}\varepsilon}{e}
-\pairing{\nabla_{\lie{a'}{a}}\varepsilon}{e} \\
&\qquad = \rho(a')\pairing{\nabla_a \varepsilon}{e}-\pairing{\nabla_{a'}\nabla_a\varepsilon}{e} \\
&\qquad = \pairing{\nabla_a\varepsilon}{\nabla_{a'}e} \\
&\qquad = \pairing{D^{\bravo{(\nabla_{a'}e)}}(\varepsilon)}{a}
\end{align*}
and
\begin{align*}
& \pairing{\big(\nabla_{a'} D^{\bravo{e}}\big)(\varepsilon)}{j(b)} \\
&\qquad = \rho(a')\pairing{D^{\bravo{e}}(\varepsilon)}{j(b)}-\pairing{D^{\bravo{e}}(\nabla_{a'}\varepsilon)}{j(b)}
-\pairing{D^{\bravo{e}}(\varepsilon)}{\lie{a'}{j(b)}} \\
&\qquad = \rho(a')\rho\big(j(b)\big)\pairing{\varepsilon}{e} - \rho\big(j(b)\big)\pairing{\nabla_{a'}\varepsilon}{e} \\
&\qquad\quad -\pairing{D^{\bravo{e}}(\varepsilon)}{-p\lie{j(b)}{a'}+j(\nabla_{a'}b)} \\
&\qquad = \rho(a')\rho\big(j(b)\big)\pairing{\varepsilon}{e} - \rho\big(j(b)\big)\rho(a')\pairing{\varepsilon}{e}
+\rho\big(j(b)\big)\pairing{\varepsilon}{\nabla_{a'}e} \\
&\qquad\quad +\pairing{\nabla_{p\lie{j(b)}{a'}}\varepsilon}{e}-\rho\big(j(\nabla_{a'}b)\big)\pairing{\varepsilon}{e} \\
&\qquad = \rho(\lie{a'}{j(b)})\pairing{\varepsilon}{e}+\rho\big(j(b)\big)\pairing{\varepsilon}{\nabla_{a'}e}
+\rho(p\lie{j(b)}{a'})\pairing{\varepsilon}{e} \\
&\qquad\quad -\pairing{\varepsilon}{\nabla_{p\lie{j(b)}{a'}}e}-\rho\big(j(\nabla_{a'}b)\big)\pairing{\varepsilon}{e} \\
&\qquad = \rho\big(j(b)\big)\pairing{\varepsilon}{\nabla_{a'}e}-\pairing{\varepsilon}{\nabla_{p\lie{j(b)}{a'}}e} \\
&\qquad = \pairing{D^{\bravo{(\nabla_{a'}e)}}(\varepsilon)}{j(b)}-\pairing{\big(\Theta(a',e)\big)^\dagger(\varepsilon)}{j(b)}
.\end{align*}
Therefore, we obtain
\[ \nabla_{a'}D^{\bravo{e}}=D^{\bravo{(\nabla_{a'}e)}}-\big(\Theta(a',e)\big)^\dagger ,\]
which is equivalent to Equation~\eqref{gather:tp2}.
\end{proof}

\subsubsection{The abelian category $\category{A}$}

It is a classical  result that the space $\sections{A}$ of smooth sections of a Lie algebroid $A$ over a smooth manifold $M$
is a Lie-Rinehart algebra over the commutative ring $C^\infty(M)$~\cite{Huebschmann,Rinehart}.
We denote the abelian category of modules over this Lie-Rinehart algebra
by the symbol $\category{A}$.
Alternatively,  $\category{A}$ can be seen as
the  category of left modules over the universal enveloping algebra
$\enveloppante{A}$ \cite{quantumgroupoid:Xu} of the Lie algebroid $A$.
In particular, the space of smooth sections of an $A$-module,
i.e.\ a vector bundle over $M$ endowed with an $A$-action,
is an object in $\category{A}$. By $D^+(\category{A})$, we denote
the bounded below derived category of $\category{A}$, which
is a symmetric monoidal category \cite{Joyal_Street}.
The interchange isomorphism $\tau:X\varotimes Y\to Y\varotimes X$
of a pair of objects $X$ and $Y$ of $D^+ (\category{A})$
is given by
\begin{equation}\label{braiding}
\tau(x\varotimes y)=(-1)^{\abs{x}\abs{y}}y\varotimes x
.\end{equation}

\subsubsection{Extension class of the jet sequence}

A short exact sequence of $A$-modules
\begin{equation}\label{Eqt:exactAmodulesequence}
\xymatrix{0 \ar[r] & P \ar[r]^\alpha & Q \ar[r]^\beta & R \ar[r] & 0}
\end{equation}
determines an extension class in the group $\text{Ext}^1_\category{A} (R, P)$,
which is naturally isomorphic to the Lie algebroid cohomology group
$H^1(A, R^*\otimes P)$ \cite{Huebschmann}.

Indeed, given a homomorphism of vector bundles $s:R\to Q$
such that $\beta\rond s=\id_R$, we have
\[ s_\sharp(\nabla_a r)-\nabla_a \big(s_\sharp(r)\big)\in\ker\beta,
\qquad \forall a\in\sections{A}, r\in\sections{R} \]
so that the equation
\begin{equation}\label{threestars}
s_\sharp(\nabla_a r)-\nabla_a \big(s_\sharp(r)\big)=\alpha_\sharp\big(\xi_s(a)\cdot r\big)
\end{equation}
defines a vector bundle map $\xi_s:A\to\Hom(R,P).$
Rewriting \eqref{threestars} as
\begin{equation}\label{fourstars} (\id_{R^*}\otimes\alpha)\rond\xi_s=\partial^A s \end{equation}
and recalling that $\alpha:P\to Q$ is a morphism of $A$-modules,
we immediately see that
\[ (\id_{R^*}\otimes\alpha)\rond\partial^A\xi_s
=\partial^A\big((\id_{R^*}\otimes\alpha)\rond\xi_s\big)
=\partial^A(\partial^A s)=0 .\]
Therefore $\partial^A\xi_s=0$, i.e.\ $\xi_s$ is a 1-cocycle for the Lie algebroid $A$
with values in the $A$-module $R^*\otimes P$.
It follows from Equation~\eqref{fourstars} that
the cohomology class $[\xi_s]\in H^1(A, R^*\otimes P)$
of the 1-cocycle $\xi_s$ defined by~\eqref{threestars}
is independent of the choice of the section $s:R\to Q$.
In fact, $[\xi_s]$ is the extension class
in $\Ext^1_\category{A}(R, P)\isomorphism H^1(A, R^*\otimes P)$
of the short exact sequence of $A$-modules \eqref{Eqt:exactAmodulesequence}.

\begin{proposition}
Given a Lie pair $(L,A)$ and an $A$-module $E$,
let $\nabla$ denote the $L$-connection on $E$ determined by
a section $s:E\to\jet{{\LoverA}}{E}$ of the short exact sequence \eqref{exact2}.
When considered as sections of $A^*\otimes A^\perp\otimes\End E$,
the bundle maps $\xi_s:A\to\Hom(E,A^\perp\otimes E)$
and $R^\nabla_E:A\otimes (L/A)\to\End E$ (respectively defined by~\eqref{threestars}
and~\eqref{Eqt:defnatyahE}) are one and the same.
\end{proposition}

\begin{proof}
Define $\breve{d}^\nabla:\sections{E}\to\sections{A^\perp\otimes E}$ by
\begin{equation}\label{shoe}
\breve{f}_\sharp(\breve{d}^\nabla e)=\bravo{e}-s_\sharp(e), \ \ \forall e\in
\Gamma (E).
\end{equation}
Then, for all $b\in B$, we have
$\pairing{\breve{d}^\nabla e}{j(b)}=\nabla_{j(b)}e$.

For all $a\in\sections{A}$, and $e\in\sections{E}$, we have
\begin{align*}
& \breve{f}_\sharp\big(\xi_s(a)\cdot e\big) \\
&\qquad = s_\sharp(\nabla_a e)-\nabla_a (s_\sharp e) && \text{by~\eqref{threestars}}, \\
&\qquad = \big(\bravo{(\nabla_a e)}-\breve{f}_\sharp\breve{d}^{\nabla}(\nabla_a e)\big)
-\big(\nabla_a(\bravo{e})-\nabla_a \breve{f}_\sharp (\breve{d}^\nabla e)\big)
&& \text{by~\eqref{shoe}}, \\
&\qquad = \breve{f}_\sharp \big(\Theta(a,e)
+\nabla_a (\breve{d}^\nabla e)-\breve{d}^{\nabla} (\nabla_a e) \big)
&& \text{by~\eqref{gather:tp1} and \eqref{gather:tp2}}.
\end{align*}
Hence, for all $b\in\sections{B}$, we get
\begin{align*}
& \pairing{\xi_s(a)\cdot e}{j(b)} \\
&\qquad = \pairing{\Theta(a,e)+\nabla_a (\breve{d}^\nabla e)-\breve{d}^{\nabla} (\nabla_a e)}{j(b)} \\
&\qquad = \nabla_{p\lie{j(b)}{a}}e+\nabla_a\pairing{\breve{d}^\nabla e}{j(b)}
-\pairing{\breve{d}^\nabla e}{j(\nabla_a b)} -\pairing{\breve{d}^\nabla (\nabla_a e)}{j(b)} \\
&\qquad = \nabla_{p\lie{j(b)}{a}}e+\nabla_a\nabla_{j(b)}e
-\nabla_{j(\nabla_a b)}e-\nabla_{j(b)}\nabla_a e \\
&\qquad = \nabla_a\nabla_{j(b)}e-\nabla_{j(b)}\nabla_a e
-\nabla_{j(\nabla_a b)+p\lie{a}{j(b)}}e \\
&\qquad = \nabla_a\nabla_{j(b)}e-\nabla_{j(b)}\nabla_a e
-\nabla_{\lie{a}{j(b)}}e \\
&\qquad = R^\nabla\big(a,j(b)\big)e \\
&\qquad = R^\nabla_E(a;b)\cdot e
.\end{align*}
This proves that $\xi_s=R^\nabla_E$.
\end{proof}

\begin{corollary}
Let $(L,A)$ be a Lie pair, and $E$ an $A$-module.
\begin{enumerate}
\item A section $s:E\to\jet{{\LoverA}}{E}$ of the short exact sequence \eqref{exact2}
is a morphism of $A$-modules if and only if the $L$-connection it induces on $E$
is compatible with the $A$-action on $E$.
\item The short exact sequence of $A$-modules \eqref{exact2} splits
if and only if the Atiyah class $\atiyahclass_E$ vanishes.
\end{enumerate}
\end{corollary}

\begin{theorem}
Let $(L,A)$ be a Lie pair, and $E$ an $A$-module.
The natural isomorphism
\[ \Ext^1_\category{A}(E,A^\perp\otimes E) \xto{\cong} H^1(A, A^\perp\otimes\End E) \]
maps the extension class of the short exact sequence of $A$-modules \eqref{exact2}
to the Atiyah class of $E$.
\end{theorem}

We refer the reader to~\cite[Lemma~8.2.4]{Calaque} for a related result
regarding the Atiyah class of dDG algebras, which correspond to the matched pair case
as pointed out in Remark~\ref{rmk:dDG}.

\section{\leibnizinfinity\ algebras}
\label{theorems}

In this section, we will explore the rich algebraic structures underlying
the Atiyah class of a Lie pair.
As we will see in the subsequent discussion,
the adequate framework is the notion of \leibnizinfinity\ algebras.
Loday's \leibnizinfinity\ algebras \cite{Loday} are a natural generalization of
Stasheff's $L_\infty$ algebras \cite{Lada_Stasheff,Lada_Markl},
where the (skew-)symmetry requirement is dropped.

Throughout this section, we implicitly identify objects of $\category{A}$
to  complexes in $\category{A}$ concentrated in degree 0.
Moreover, we make frequent use of the shifting functor:
the shift $V[k]$ of a graded vector space $V=\bigoplus_n V_n$ is determined
by the rule $\big(V[k]\big)_n=V_{k+n}$.

We defer most proofs to Section~\ref{Subsectin:proofsofmaintheorems}.

\subsection{$L/A[-1]$ as  a Lie algebra object}

Recall that a \emph{Lie algebra object} in a monoidal category $\category{C}$
is an object $\Lambda$ of $\category{C}$ together with a morphism
$\lambda\in\Hom_{\category{C}}(\Lambda\otimes\Lambda,\Lambda)$
such that
\begin{enumerate}
\item $\lambda\rond\tau=-\lambda$ (skew-symmetry), and
\item $\lambda\rond(\id\otimes\lambda)= \lambda\rond(\lambda\otimes\id)
+\lambda\rond(\id\otimes\lambda)\rond(\tau\otimes\id)
\qquad \text{(Jacobi identity)},$
where $\tau:\Lambda\otimes\Lambda\to\Lambda\otimes\Lambda$ is the braiding isomorphism.
\end{enumerate}

Let $(L,A)$ be a Lie pair with quotient $B=L/A$.
Note that
\[ \Hom_{D^+(\category{A})}(B\otimes B,\shiftby{B}{1}) \cong
\Hom_{D^+(\category{A})}(\shiftby{B}{-1}\otimes\shiftby{B}{-1},\shiftby{B}{-1}) .\]
Being an element of
\begin{multline*}
\Ext^1_{\category{A}}(B,B^*\otimes B) \cong
\Ext^1_{\category{A}}(B\otimes B,B) \cong
\Hom_{D^+(\category{A})}(B\otimes B,\shiftby{B}{1}) \cong \\
\Hom_{D^+(\category{A})}(\shiftby{B}{-1}\otimes\shiftby{B}{-1},\shiftby{B}{-1})
,\end{multline*}
the Atiyah class $\atiyahclass_{B}$ of the $A$-module $B$
defines a ``Lie bracket'' on $\shiftby{B}{-1}$.
If, moreover, $E$ is an $A$-module, its Atiyah class
\[ \atiyahclass_{E}\in\Ext^1_{\category{A}}(E,B^*\otimes E)
\cong\Ext^1_{\category{A}}(B\otimes E,E)\cong
\Hom_{D^+(\category{A})}(\shiftby{B}{-1}\otimes\shiftby{E}{-1},\shiftby{E}{-1}) \]
defines a ``representation'' on $\shiftby{E}{-1}$ of the ``Lie algebra'' $\shiftby{B}{-1}$.

More precisely, we have the following theorem, whose proof is deferred
to Section \ref{Subsectin:proofsofmaintheorems}.

\begin{theorem}\label{laoidc}
Let $(L,A)$ be a Lie pair with quotient $B=L/A$. Then $\shiftby{B}{-1}$
is a Lie algebra object in the derived category $D^+(\category{A})$.
Moreover, if $E$ is an $A$-module, then $\shiftby{E}{-1}$ is a module object
over the Lie algebra object $\shiftby{B}{-1}$ in the derived category $D^+(\category{A})$.
\end{theorem}

\begin{remark}
From the skew-symmetric property of a Lie algebra,
it follows that the Atiyah class $\alpha_{B}$ can  indeed
be considered as an element in $H^1(A, S^2 B^*\otimes B)$, or more precisely,
in  the  image of the map $H^1(A, S^2 B^*\otimes B)\to
 H^1(A,  B^*\otimes \End B)$
  induced by the $A$-modules morphism
$S^2 B^*\otimes B\to B^*\otimes B^* \otimes B (\cong B^*\otimes \End B)$.
\end{remark}

\begin{remark}
It is implicitly stated in~\cite{Kapranov}
(see also \cite{Willerton, Ramadoss})
that, if $X$ is a complex manifold, then $T_X[-1]$ is a Lie algebra
object in the  bounded below derived category $D^+(X)$
of  coherent sheaves on $X$.
This is simply Theorem~\ref{laoidc} in the special case when
 $L=T_X\otimes\CC$ and $A=T\zo_X$.
\end{remark}

\subsection{Jacobi identity up to homotopy}

Let $(L,A)$ be a Lie pair and $E$ an $A$-module.
The quotient $B={\LoverA}$ is naturally an $A$-module
(see Proposition~\ref{Prop:quotientmodule}).

Consider the graded vector spaces
\[ V=\bigoplus_{n=0}^\infty \sections{\wedge^n A^*\otimes B} \]
and
\[ W=\bigoplus_{n=0}^\infty \sections{\wedge^n A^*\otimes E} ,\]
and the covariant differentials
\begin{gather*}
\partialA:\Gamma\big(\wedge^{\bullet}A^*\otimes B\big)
\to\Gamma\big(\wedge^{\bullet+1}A^*\otimes B\big) \\
\partialA:\Gamma\big(\wedge^{\bullet}A^*\otimes E\big)
\to\Gamma\big(\wedge^{\bullet+1}A^*\otimes E\big)
\end{gather*}
associated to the $A$-actions on $B$ and $E$, respectively.
Choosing an $L$-connection $\nabla$ on $\bB$ extending the $A$-action,
we obtain  the bundle maps $R_2:B\otimes B\to\Hom(A,B)$
and $S_2:B\otimes E\to\Hom(A,E)$ given by
\begin{gather}
A\ni a\xmapsto{R_2(b_1,b_2)}R^\nabla_B(a;b_1)b_2 \in B , \label{r2} \\
A\ni a\xmapsto{S_2(b,e)}R^\nabla_E(a;b)e \in E , \label{s2}
\end{gather}
where $R^\nabla_B:A\otimes B\to\End B$ and $R^\nabla_E:A\otimes B\to\End E$
denote the Atiyah cocycles of $B$ and $E$, respectively.

\begin{theorem}\label{Thm:CsuperLiealgebraobject}
Up to homotopies, the complex $(V[-1],\partial^A)$ is a differential graded Lie algebra
and the complex $(W[-1],\partial^A)$ is a differential graded module over it.
The Lie algebra bracket
\[ V[-1]\varotimes V[-1]\xto{\lambda} V[-1] \]
and the representation
\[ V[-1]\varotimes W[-1]\xto{\mu} W[-1] \]
are given by
\begin{equation*}
\lambda\big( (\xi_1\otimes b_1),  (\xi_2\otimes b_2) \big)
=\minuspower{k_2} \xi_1\wedge\xi_2\wedge R_2(b_1,b_2)
\end{equation*}
and
\begin{equation*}
\mu\big( (\xi_1\otimes b),  (\xi_2\otimes e) \big)
=\minuspower{k_2} \xi_1\wedge\xi_2\wedge S_2(b,e)
,\end{equation*}
where $\xi_1\in\sections{\wedge^{k_1}A^*}$,
 $\xi_2\in\sections{\wedge^{k_2}A^*}$,
$b_1,b_2,b\in\sections{B}$, and $e\in\sections{E}$.
\end{theorem}
Consequently, the cohomology $\bigoplus_{i\geq 1}H^{i-1}(A,E)=H^\bullet(W[-1],\partial^A)$
is a module over the (graded) Lie algebra $\bigoplus_{i\geq 1}H^{i-1}(A,B)=H^\bullet(V[-1],\partial^A)$.


In Section \ref{Subsection:homotopymaintheorems}, we will
describe a result which keeps track of higher homotopies.

\subsection{\leibnizinfinityone\ algebras}

Recall that a graded Leibniz algebra  is a $\ZZ$-graded vector space
$V=\bigoplus_{k\in\ZZ}V_k$
equipped with a bilinear bracket $V\otimes V\xto{\ba{-}{-}}V$ of degree 0
satisfying the graded Leibniz rule
\begin{equation*}
\ba{x}{\ba{y}{z}}=\ba{\ba{x}{y}}{z}
+(-1)^{\abs{x}\abs{y}}\ba{y}{\ba{x}{z}},
\end{equation*}
for all homogeneous elements $x,y,z\in V$.

If, moreover, $V$ is endowed with a differential $\delta$ of degree $1$ satisfying
\[ \delta\ba{x}{y}=\ba{\delta x}{y}+\minuspower{\abs{x}+1}\ba{x}{\delta y} \]
for all homogeneous elements $x,y\in V$,
then we say that $(V,\ba{-}{-},\delta)$ is a
 differential graded Leibniz algebra.

\begin{definition}
A \leibnizinfinityone\ algebra is a $\ZZ$-graded vector space
$V=\bigoplus_{n\in\ZZ}V_n$ endowed with
a sequence $(\lambda_k)_{k=1}^\infty$ of multilinear maps
$\lambda_k: \otimes^k V\to V$ of degree $1$
satisfying the identity

\begin{multline}\label{Eqt:generalJacobi}
\sum_{1\le j\le k\le n}\sum_{\sigma\in\shuffle{k-j}{j-1}}
\Koszul{\sigma;v_1,\cdots,v_{k-1}}
(-1)^{\abs{v_{\sigma(1)}}+\abs{v_{\sigma(2)}}+\cdots+\abs{v_{\sigma(k-j)}}} \\
\lambda_{n-j+1}\big(v_{\sigma(1)},\cdots,v_{\sigma(k-j)},
\lambda_j(v_{\sigma(k+1-j)},\cdots,v_{\sigma(k-1)},v_k),
v_{k+1},\cdots, v_{n}\big)=0
\end{multline}

for each $n\in\NN$ and for any homogeneous vectors $v_1,v_2,\dots,v_n\in V$.
Here $\shuffle{k-j}{j-1}$ denotes the set of
$(k-j,j-1)$-shuffles\footnote{A $(k-j, j-1)$-shuffle is a permutation $\sigma$
of the set $\{1,2,\cdots,k-1\}$ such that $\sigma(1)\le\sigma(2)\le\cdots\le\sigma(k-j)$
and $\sigma(k-j+1)\le\sigma(k-j+2)\le\cdots\le\sigma(k-1)$.},
and $\Koszul{\sigma; v_1, \cdots, v_{k-1}}$ denotes
the Koszul sign\footnote{The Koszul sign of a permutation $\sigma$
of the (homogeneous) vectors $v_1,v_2,\dots,v_n$ is determined by
the relation $v_{\sigma(1)}\odot v_{\sigma(2)}\odot\cdots\odot v_{\sigma(n)}
= \Koszul{\sigma; v_1, \cdots, v_n} \cdot v_1\odot v_2\odot\cdots\odot v_n$.}
of the permutation $\sigma$ of the (homogeneous) vectors $v_1,v_2,\dots,v_{k-1}$.
\end{definition}

\begin{remark}
If all $\lambda_k$ are zero except for $\lambda_1$,
$(V,\lambda_1)$ is simply a cochain complex.
If $\lambda_k=0$ ($k\geq 3$), then
$(\shiftby{V}{-1},\ba{-}{-}, d)$ is a graded
differential Leibniz algebra, where
$\ba{x}{y}=\minuspower{\abs{x}}\lambda_2(x,y)$, and $d=\lambda_1$.
\end{remark}

\begin{remark}
A graded vector space $V$ is a \leibnizinfinityone\ algebra if and only if the shifted graded vector space $V[-1]$
is a \leibnizinfinity\ algebra in the sense of Loday
\cite{Poncin,Uchino}.
Working with \leibnizinfinityone\ algebras rather than \leibnizinfinity\ algebras
is convenient,
 as all maps in the sequence $(\lambda_k)_{k=1}^n$ have the same degree
in this setting.
\end{remark}

\begin{definition}
A module over a \leibnizinfinityone\ algebra $V$ is a $\ZZ$-graded vector space
$W=\bigoplus_{n\in\ZZ}W_n$ together with a sequence $(\mu_k)_{k=1}^\infty$
of multilinear maps \[ \mu_k: (\otimes^{k-1} V)\otimes W\to W \] of degree $1$
satisfying the identity:
\begin{multline*}
\sum_{1\le j\le k\le n-1}\sum_{\sigma\in\shuffle{k-j}{j-1}}
\Koszul{\sigma;v_1,\cdots,v_{k-1}}
(-1)^{\abs{v_{\sigma(1)}}+\abs{v_{\sigma(2)}}+\cdots+\abs{v_{\sigma(k-j)}}} \\
\mu_{n-j+1}\big(v_{\sigma(1)},\cdots,v_{\sigma(k-j)},
\lambda_j(v_{\sigma(k+1-j)},\cdots,v_{\sigma(k-1)},v_k),
v_{k+1},\cdots, v_{n-1},w\big) \\
+ \sum_{1\le j\le n}\sum_{\sigma\in\shuffle{n-j}{j-1}}
\Koszul{\sigma;v_1,\cdots,v_{n-1}}
(-1)^{\abs{v_{\sigma(1)}}+\abs{v_{\sigma(2)}}+\cdots+\abs{v_{\sigma(n-j)}}} \\
\mu_{n-j+1}\big(v_{\sigma(1)},\cdots,v_{\sigma(n-j)},
\mu_j(v_{\sigma(n+1-j)},\cdots,v_{\sigma(n-1)},w)\big) =0,
\end{multline*}
for each $n\in\NN$ and any homogeneous vectors $v_1,v_2,\dots,v_{n-1}\in V$
and $w\in W$.
\end{definition}

\begin{remark}
A graded vector space $W$ is a module over a \leibnizinfinityone\ algebra $V$ if and only if
$V\oplus W$ is a  \leibnizinfinityone\ algebra such that
$V$ is a \leibnizinfinityone\ subalgebra \cite{Lada}.
\end{remark}

The proof of the next proposition is a direct verification, which we omit.
\begin{proposition}
If $\big(V,(\lambda_k)_{k=1}^\infty\big)$ is a \leibnizinfinityone\ algebra,
then $(V, \lambda_1)$ is a cochain complex and its cohomology $H\graded(V)[-1]$
is a graded Leibniz algebra with bracket $H(\lambda_2)$, the image of $\lambda_2$
(seen as a chain map) under the cohomology functor.
Moreover, if $\big(W,(\mu_k)_{k=1}^\infty\big)$ is a module over $\big(V,(\lambda_k)_{k=1}^\infty\big)$,
then $(W,\mu_1)$ is a cochain complex and $H(\mu_2)$ is a representation of $H\graded(V)[-1]$ on the
cohomology $H\graded(W)[-1]$ of $(W,\mu_1)$.
\end{proposition}

\subsection{Main theorem}\label{Subsection:homotopymaintheorems}

Unless we state otherwise, we assume throughout this section that
$(L,A)$ is a Lie pair and $E$ is an $A$-module.
The quotient $B={\LoverA}$ is naturally an $A$-module
(see Proposition~\ref{Prop:quotientmodule}).
We use the symbol $\partial^A$ to denote the covariant differential
\[ \partialA:\Gamma\big(\wedge^{\bullet}A^*\otimes(\otimes^{\star}B^*)\otimes E\big)
\to\Gamma\big(\wedge^{\bullet+1}A^* \otimes(\otimes^{\star}B^*)\otimes E\big) \]
associated to the $A$-action on $(\otimes^{\star}B^*)\otimes E$.
In particular, for any bundle map
$\mu:(\wedge^k A)\otimes(\otimes^l B)\to B$, we have
\begin{multline*}
\big(\partial^A\mu\big)(a_0\wedge\cdots\wedge a_k;b_1\otimes\cdots\otimes b_l) = \\
\sum_{i=0}^k (-1)^i \left\{
\nabla_{a_i} \big( \mu(a_{\widehat{i}};b_1\otimes\cdots\otimes b_l)\big)
-\mu\big(a_{\widehat{i}};\nabla_{a_i} (b_1\otimes\cdots\otimes b_l) \big) \right\} \\
+ \sum_{i<j} (-1)^{i+j} \mu(\lie{a_i}{a_j}\wedge
a_{\widehat{i,j}};b_1\otimes\cdots\otimes b_l)
,\end{multline*}
where $a_{\widehat{i}}$ stands for
$a_0\wedge\cdots\wedge\widehat{a_i}\wedge\cdots\wedge a_k$ and
$a_{\widehat{i,j}}$ for
$a_0\wedge\cdots\wedge\widehat{a_i}\wedge\cdots\wedge\widehat{a_j}\wedge
\cdots\wedge a_k$, and $\nabla_{a_i}(b_1\otimes\cdots\otimes b_l)$ for
$\sum_{j=1}^l b_1\otimes\cdots\otimes\nabla_{a_i}b_j\otimes\cdots\otimes b_l$.

\subsubsection{The operator $\partial^\nabla$}
\label{section:341}
Now choose an extension of the $A$-action on $E$ to an $L$-connection $\nabla$ on $E$,
an extension of the $A$-action on $B$ to an $L$-connection $\nabla$ on $B$,
and a splitting of the short exact sequence of vector bundles
\begin{equation}\label{ABL}
\xymatrix{0 \ar[r] & A \ar[r]^i & L \ar[r]^q & B \ar[r] & 0 }
,\end{equation}
i.e.\ a pair of maps $j:B\to L$ and $p:L\to A$ such that
$q\rond j=\id_B$, $p\rond i=\id_A$ and $i\rond p+j\rond q=\id_L$:
\[ \xymatrix{0 \ar@<.5ex>[r] & A \ar@<.5ex>[r]^i \ar@<.5ex>[l] & L \ar@<.5ex>[r]^q \ar@<.5ex>[l]^p
& B \ar@<.5ex>[r] \ar@<.5ex>[l]^j & 0 \ar@<.5ex>[l] } .\]
This splitting determines a  map
\[ \sections{B}\times\sections{A}\to\sections{A}: (b,a) \mapsto p\lie{j(b)}{i(a)} ,\]
which we will denote by $\anadelta$ since it satisfies the relations
\[ \anadelta_{f b}a=f \anadelta_{b}a
\qquad\text{and}\qquad
\anadelta_{b}(fa)=\pairing{\rho^* df}{j(b)} a+f\anadelta_b a ,\]
for all $f\in\cinf{M}$, $b\in\sections{B}$, and $a\in\sections{A}$.
In some sense, $B$ ``acts'' on $A$.
\newline
Identifying sections of $\wedge^\bullet A^*\otimes(\otimes^\star B^*)\otimes E$
with bundle maps $\wedge^\bullet A\otimes(\otimes^\star B)\to E$, we define
a differential operator
\begin{equation}
\label{eq:cat}
 \partial^\nabla:\sections{\wedge^k A^* \otimes (\otimes^l B^*) \otimes E}
\to \sections{\wedge^k A^* \otimes (\otimes^{l+1} B^*) \otimes E}
\end{equation}
by \[ (-1)^k \interior{b_0}(\partial^\nabla\omega)=\nabla_{j(b_0)}\omega \] or, more precisely,
\begin{multline*}
\minuspower{k} \big(\partial^\nabla\omega\big)(a_1,\cdots,a_k;b_0,\cdots,b_l)
= \nabla_{j(b_0)}\big(\omega(a_1,\cdots,a_k;b_1,\cdots,b_l)\big) \\
-\omega(\anadelta_{b_0}a_1,\cdots,a_k;b_1,\cdots,b_l) -\cdots
-\omega(a_1,\cdots,\anadelta_{b_0}a_k;b_1,\cdots,b_l) \\
-\omega(a_1,\cdots,a_k;\nabla_{j(b_0)}b_1,\cdots,b_l) -\cdots
-\omega(a_1,\cdots,a_k;b_1,\cdots,\nabla_{j(b_0)}b_l)
,\end{multline*}
where $a_1,\cdots,a_k\in\sections{A}$, $b_0,\cdots,b_l\in\sections{B}$,
and $\omega:\wedge^k A\otimes(\otimes^l B)\to E$.

Note that $\partial^\nabla$ depends on the choice of the $L$-connections extending the $A$-actions
and the splitting $j:B\to L$ of the short exact sequence \eqref{ABL}, while $\partial^A$ does not.

The chosen splitting of \eqref{ABL} does also determine three vector bundle maps
\begin{equation}
\label{eq:snake}
 \alpha:\wedge^2 B\to A, \qquad \beta:\wedge^2 B\to B, \quad \text{and} \quad
\Omega:\wedge^2 B\to \End B
\end{equation}
 given by
\begin{gather}
\alpha(b_1,b_2)=p\lie{j(b_1)}{j(b_2)} \label{eq:snake1} \\
\beta(b_1,b_2)=\nabla_{j(b_1)}b_2-\nabla_{j(b_2)}b_1-q\lie{j(b_1)}{j(b_2)} \label{eq:snake2}\\
\intertext{and}
\Omega(b_1,b_2)=\nabla_{j(b_1)}\nabla_{j(b_2)}-\nabla_{j(b_2)}\nabla_{j(b_1)}
-\nabla_{\lie{j(b_1)}{j(b_2)}}
. \label{eq:snake3}
\end{gather}

\begin{proposition}\label{Lemma1}
For any $a\in\sections{A}$ and $b_1,b_2\in\sections{B}$, we have
\[ R^\nabla_B(a;b_1)b_2-R^\nabla_B(a;b_2)b_1=\big(\nabla_a\beta\big)(b_1,b_2) \]
or, equivalently,
\[ R_2(b_1,b_2)-R_2(b_2,b_1)=\big(\partial^A \beta\big)(b_1,b_2) .\]
\end{proposition}

\begin{proof}
For convenience, we set $\widetilde{b}=j(b)$, $\forall b\in \Gamma( B)$.
Hence $\lie{a}{\widetilde{b}}=-\anadelta_b a+\widetilde{\nabla_a b}$ and
\[ \lie{\widetilde{b_1}}{\widetilde{b_2}}=\alpha(b_1,b_2)
+\widetilde{\nabla_{\widetilde{b_1}} b_2}-\widetilde{\nabla_{\widetilde{b_2}} b_1}
-\widetilde{\beta(b_1,b_2)} .\]
A straightforward computation yields the equality
\begin{multline*} q\big(\lie{\lie{a}{\widetilde{b_1}}}{\widetilde{b_2}}
+\lie{\lie{\widetilde{b_1}}{\widetilde{b_2}}}{a}
+\lie{\lie{\widetilde{b_2}}{a}}{\widetilde{b_1}}\big) \\
=R^\nabla_B(a;b_2)b_1-R^\nabla_B(a;b_1)b_2
+\nabla_a\big(\beta(b_1,b_2)\big)-\beta(\nabla_a b_1,b_2)-\beta(b_1,\nabla_a b_2)
.\end{multline*}
The result follows from the Jacobi identity of the Lie algebroid $L$.
\end{proof}

Note that, since $R^\nabla_B$ is (by its very definition) independent of the choice of the splitting,
Proposition~\ref{Lemma1} asserts that, unlike $\beta$, $\partial^A\beta$ does not depend on the choice of splitting.

\subsubsection{The maps $R_n$}

Recall the bundle map $R_2:B\otimes B\to\Hom(A,B)$
associated to the Atiyah cocycle of $B$ given by \eqref{r2}.
Since $B$ is an $A$-module, we can substitute $B$ for $E$
in Equation \eqref{eq:cat}
and define a sequence of bundle maps \[ R_n:\otimes^n B\to\Hom(A,B) \]
inductively by the relation
\begin{equation}
\label{eq:dog}
 R_{n+1}=\partial^\nabla R_n,  \qquad \text{for } n\geq 2 .
\end{equation}
Hence, we have
\[ R_{n+1}(b_0\otimes b_1\otimes\cdots\otimes b_n)
=R_n\big(\nabla_{j(b_0)}(b_1\otimes\cdots\otimes b_n)\big)
-\nabla_{j(b_0)}\big(R_n(b_1\otimes\cdots\otimes b_n)\big) .\]

\begin{example}\label{Ex:ABmathchedpairtrivialexample}
Let $L=A\bowtie B$ be a matched pair of Lie algebras.
Any bilinear map $\gamma:B\otimes B\to B$ determines an $L$-connection $\nabla$ on $B$
extending its $A$-module structure (and conversely): $\nabla_{b_1}b_2=\gamma(b_1,b_2)$.
Taking $\gamma=0$, the Atiyah cocycle reads
\[ \atiyahcocycle_{B}(a;b_1)b_2=\nabla_{\Delta_{b_1}a}b_2. \]
Hence
 \[ R_n(b_1,b_2,b_3,\cdots,b_n)=
\nabla_{\Delta_{b_{n-1}}\Delta_{b_{n-2}}\cdots
\Delta_{b_{1}}(\argument)}b_n .\]
\end{example}

\subsubsection{\leibnizinfinityone\ algebra (and modules) arising from a Lie pair}

Consider the sequence of $k$-ary operations
$\lambda_k:\otimes^kV\to V$ ($k\in\NN$)
on the graded vector space
\[ V=\bigoplus_{n=0}^\infty \sections{\wedge^n A^*\otimes B} \]
defined by $\lambda_1=\partial^A$  and, for $k\geq 2$,
\begin{equation}\label{eqn:temp1}
\lambda_k(\xi_1\otimes b_1,\cdots,\xi_k\otimes b_k)
=\minuspower{\abs{\xi_1}+\cdots+\abs{\xi_k}}
\xi_1\wedge\cdots\wedge\xi_k\wedge R_k(b_1,\cdots, b_k)
,\end{equation}
where $b_1,\dots,b_k\in\sections{B}$ and $\xi_1,\dots,\xi_k$
are homogeneous elements in $\sections{\wedge^\bullet A^*}$.

We are now ready to state the main result of the paper.

\begin{theorem}\label{Thm:CisSHLeiniz}
When endowed with the sequence of multibrackets $(\lambda_k)_{k\in\NN}$ defined above,
the graded vector space $V=\bigoplus_{n=0}^\infty \sections{\wedge^n A^*\otimes B}$
becomes a \leibnizinfinityone\ algebra.
\end{theorem}

Similarly, we can introduce the bundle map $S_2:B\otimes E\to\Hom(A,E)$ given by
\[ A\ni a\xmapsto{S_2(b;e)}R^\nabla_E(a;b)\cdot e\in E,\]
where $R^\nabla_E:A\otimes B\to\End E$ denotes the Atiyah cocycle of the $A$-module $E$,
and then define a sequence of bundle maps \[ S_n:(\otimes^{n-1} B)\otimes E\to\Hom(A,E) \]
inductively by the relation \[ S_{n+1}=\partial^\nabla S_n,  \qquad \text{for } n\geq 2 .\]

Consider the graded vector space $W=\bigoplus_{n=0}^\infty \sections{\wedge^n A^*\otimes E}$
and the sequence of $k$-ary brackets
$\mu_k:(\otimes^{k-1}V) \otimes W \to W$ ($k\in\NN$)
defined by $\mu_1=\partial^A$  and, for $k\geq 2$,
\begin{equation*}
\mu_k(\xi_1\otimes b_1,\cdots,\xi_{k-1}\otimes b_{k-1};\xi_k\otimes e)
=\minuspower{\abs{\xi_1}+\cdots+\abs{\xi_{k}}}
\xi_1\wedge\cdots\wedge\xi_k\wedge S_k(b_1,\cdots,b_{k-1};e)
,\end{equation*}
where $b_1,\dots,b_{k-1}\in\sections{B}$, $e\in\sections{E}$, and $\xi_1,\dots,\xi_k$
are homogeneous elements of $\sections{\wedge^\bullet A^*}$.

\begin{theorem}\label{thm:module}
When endowed with the sequence of multibrackets $(\mu_k)_{k\in\NN}$ defined above,
the graded vector space $W=\bigoplus_{n=0}^\infty \sections{\wedge^n A^*\otimes E}$
becomes a \leibnizinfinityone\ module over the \leibnizinfinityone\ algebra $\big(V,(\lambda_k)_{k\in\NN}\big)$.
\end{theorem}

\begin{example}
A Lie bialgebra $(\mfg,\mfg^*)$ is a matched pair of Lie algebras.
Therefore, it induces two Lie pairs: $(\mfg\bowtie\mfg^*,\mfg)$ and $(\mfg\bowtie\mfg^*,\mfg^*)$.
It follows from Example~\ref{Ex:ABmathchedpairtrivialexample} and Theorem~\ref{Thm:CisSHLeiniz}
that both $\bigoplus_{n\geq 0}\wedge^n\mfg^*\otimes\mfg^*$
and $\bigoplus_{n\geq 0}\wedge^n\mfg\otimes\mfg$ are \leibnizinfinityone\ algebras.
\end{example}

Let $A$ be a Lie algebroid over a manifold $M$.
By an $A$-algebra, we mean a bundle (of finite or infinite rank)
of associative algebras $\Cp$ over $M$, which is an $A$-module,
and on which $\Gamma (A)$ acts by derivations.
For a commutative $A$-algebra $\Cp$, the sequence of maps $(\lambda_k)_{k\in\NN}$
extends,  in a natural way,  to the graded vector  space
$\bigoplus_{n=0}^\infty\sections{\wedge^n A^*\otimes B\otimes \Cp}$.
Similarly, the sequence of maps $(\mu_k)_{k\in\NN}$ extends  to the graded space
$\bigoplus_{n=0}^\infty\sections{\wedge^n A^*\otimes E\otimes \Cp}$.

\begin{theorem}
\label{thm:mainC}
Let $(L,A)$ be a Lie pair with quotient $B$, and let $\Cp$ be a commutative $A$-algebra.
When endowed with the sequence of multibrackets $(\lambda_k)_{k\in\NN}$,
the graded vector space
$\sections{\wedge^\bullet A^*\otimes B\otimes \Cp}$
becomes a Leibniz$_\infty [1]$  algebra.
Moreover, if $E$ is an $A$-module, the graded vector space
$\sections{\wedge^\bullet A^*\otimes E \otimes \Cp}[-1]$
becomes a  Leibniz$_\infty [1]$ module over the
Leibniz$_\infty [1]$ algebra
$\sections{\wedge^\bullet A^*\otimes B \otimes \Cp}[-1]$.
\end{theorem}

As an immediate  consequence, we have the following

\begin{corollary}\label{cor:Lie}
Under the same hypothesis as in Theorem \ref{thm:mainC}, $\bigoplus_{i\geq 1}H^{i-1} (A,B\otimes\Cp)$
is a graded Lie algebra,
 and $\bigoplus_{i\geq 1}H^{i-1}(A, E \otimes \Cp)$ a module over it.
\end{corollary}

\begin{example}
Let $\mfg$ be a Lie subalgebra of a Lie algebra $\mfd$ as in
Example~\ref{Ex:dhliealgebras}. Assume that $\Cp$ is a commutative
$\mfg$-algebra.
\newline
Every linear map $\boldsymbol{L}: \mfd\to\End(\mfd/\mfg)$ that extends
 the $\mfg$-module structure $\mfg\to\End(\mfd/\mfg)$ induces a
2-ary bracket on $ \wedge^{\bullet-1}\mfg^* \otimes \mfd/\mfg\otimes \Cp$:
\begin{equation}
\label{2ary}
\left[\classe{\xi_1\otimes b_1\otimes c_1},\classe{\xi_2\otimes b_2 \otimes c_2}\right]
=\minuspower{\abs{\xi_2}}\classe{\xi_1\wedge\xi_2\otimes
\big(\partial^{\mfg}\boldsymbol{L}\big)(\argument;b_1)\cdot b_2} \otimes c_1c_2,
\end{equation}
which in turn induces a (graded) Lie algebra bracket
on the Chevalley-Eilenberg cohomology
 $\bigoplus H^{\bullet-1}(\mfg, \mfd/\mfg \otimes\Cp)$.
Here $\xi_i\otimes b_i\otimes c_i$ ($i=1,2$) are cocycles with
$\xi_1,\xi_2 \in\wedge^\bullet\mfg^*$, $b_1,b_2\in\mfd/\mfg$ and $c_1, c_2 \in \Cp$.
\newline
Moreover, if $E$ is a $\mfg$-module, every linear map $\boldsymbol{M}:\mfd\to\End E$
that extends the $\mfg$-module structure $\mfg\to\End E$ gives rise to a
bilinear map
\[ \classe{(\xi_1\otimes b\otimes c_1)}\rhd\classe{(\xi_2\otimes e \otimes c_2)}
=\minuspower{\abs{\xi_2}}\classe{\xi_1\wedge\xi_2
\otimes \big(\partial^{\mfg}\boldsymbol{M}\big)(\argument;b)\cdot e} \otimes c_1c_2, \]
which induces a representation on
 $\bigoplus H^{\bullet-1}(\mfg, E \otimes \Cp)$
of the graded Lie algebra $\bigoplus H^{\bullet-1}(\mfg, \mfd/\mfg \otimes \Cp)$.
Here $\xi_1\otimes b\otimes c_1$ and $\xi_2\otimes e \otimes c_2$
 are cocycles with $\xi_1,\xi_2 \in\wedge^\bullet\mfg^*$, $b\in\mfd/\mfg$,  $e\in E$
and $c_1, c_2 \in \Cp$.

Take a complement  $\mathfrak{h}$ of $\mathfrak{g}$ in $\mathfrak{d}$
so that we can write $\mathfrak{d}=\mathfrak{g}\oplus\mathfrak{h}$.
Then $\mathfrak{d}/\mathfrak{g}$ can be identified with $\mathfrak{h}$,
on which  the $\mathfrak{h}$-action is given by
$a\cdot h= \prh [a, h]$. Take
$\boldsymbol{L}: \mfd\to\End\mathfrak{h}$ to be the trivial extension
of the $\mathfrak{g}$-module structure $\mathfrak{g} \to \End\mathfrak{h}$,
i.e.\ set $\boldsymbol{L}|_\mathfrak{h}=0$.
Then the 2-ary bracket in Equation \eqref{2ary} is given by
\[ [f\otimes c_1, \ g\otimes c_2]=[f,g] \otimes c_1c_2 ,\]
where  $f$ is a linear map from $\wedge^p\mathfrak{g}$ to $\mathfrak{h}$,
$g$ is  a linear map from $\wedge^q\mathfrak{g}$ to $\mathfrak{h}$,
 $c_1, c_2\in \Cp$, and
$[f,g]$ is  a linear map from $\wedge^{p+q+1}\mathfrak{g}$
to $\mathfrak{h}$  given by
\begin{multline*} [f,g](a_0,a_1,\cdots,a_{p+q})= \\
-\sum_{\sigma\in\mathfrak{S}_{1,p,q}}\sgn(\sigma)
\prh\lie{\prg\lie{a_{\sigma(0)}}{f(a_{\sigma(1)},\cdots,a_{\sigma(p)})}}
{g(a_{\sigma(p+1)},\cdots,a_{\sigma(p+q)})} .\end{multline*}
Here $\mathfrak{S}_{1,p,q}$ is the set of all permutations
$\sigma$ of $\{0,1,\cdots,p+q\}$ satisfying
$\sigma(1)<\cdots<\sigma(p)$ and $\sigma(p+1)<\cdots<\sigma(p+q)$.
\end{example}

\begin{remark}
It is natural to ask  how  the  \leibnizinfinityone\ algebra structure
obtained in Theorem \ref{Thm:CisSHLeiniz} and the \leibnizinfinityone\
module structure in Theorem \ref{thm:module} depend on
the choice of connections and the splitting data.
This question will be investigated somewhere else.
\end{remark}

\subsubsection{$\Linf$ rather than \leibnizinfinity}

\begin{theorem}
\label{Thm:conditionforCtobeSHLie}
Let $(A,B)$ be a matched pair of Lie algebroids with direct sum $L=A\bowtie B$.
Assume there exists a flat torsion free $B$-connection on $B$.
Then the maps $R_n$ defined as in Equation \eqref{eq:dog}
 are totally symmetric, the multibrackets
$\lambda_k:\otimes^k V\to V$ ($k\in\NN$) defined as
in Equation \eqref{eqn:temp1} on the graded vector space
$V=\bigoplus_{n=0}^{\infty}\sections{\wedge^n A^*\otimes B}$ are
 graded symmetric,
and $V[-1]$ is actually an $\Linf$ algebra.
\end{theorem}

The following example is due to Camille Laurent-Gengoux.

\begin{example}
The general Lie algebra $\mathfrak{gl}_n(\CC)$ decomposes as the direct sum of the unitary Lie algebra $\mathfrak{u}_n$ and the Lie algebra $\mathfrak{t}_n$ of upper triangular matrices with real diagonal coefficients.
Both  $\mathfrak{u}_n$ and $\mathfrak{t}_n$ are isotropic with respect
to the natural nondegenerate ad-invariant inner product
\[ X\otimes Y\mapsto \im\big(\tr(XY)\big) \] on
$\mathfrak{gl}_n(\CC)$.
Hence $(\mathfrak{u}_n,\mathfrak{t}_n)$ is a matched pair of Lie algebras as well as a Lie bialgebra.
Matrix multiplication being associative, setting $\nabla_X Y=XY$ for any $X,Y\in\mathfrak{t}_n$ defines
a flat torsion free $\mathfrak{t}_n$-connection on $\mathfrak{t}_n$.
It follows from Theorem~\ref{Thm:conditionforCtobeSHLie} that
$\sections{\wedge^\bullet \mathfrak{u}_n^*\otimes \mathfrak{t}_n}[-1]
\isomorphism\sections{\wedge^\bullet\mathfrak{t}_n\otimes\mathfrak{t}_n}[-1]$ is an
$L_{\infty}$ algebra.
\end{example}

As an application of Theorem  \ref{Thm:conditionforCtobeSHLie},
consider a  K\"ahler manifold $X$.
The complexification $\nabla^\CC$ of its Levi-Civita connection is a $T_X\otimes\CC$-connection on  $T_X\otimes\CC$.
Set $A=T\zo_X$ and $B=T\oz_X$.
Then $(A,B)$ is a matched pair of Lie algebroids, whose direct sum $A\bowtie B$ is isomorphic,
as a Lie algebroid, to $T_X\otimes\CC$.
It is easy to see that $\nabla^\CC$ induces a flat torsion free $B$-connection on $B$.
In this context, the tensors $R_n\in\OO^{0,1}\big(\Hom(\otimes^n \pT,\pT)\big)$
are the curvature $R_2\in\OO^{1,1}(\End \pT)\isomorphism\OO^{0,2}\big(\Hom(\pT\otimes \pT,\pT)\big)$
and its higher covariant derivatives: $R_{i+1}=\partial^\nabla R_i$.
Applying Theorem~\ref{Thm:conditionforCtobeSHLie}, we recover a result of Kapranov~\cite{Kapranov}:

\begin{corollary}[Kapranov]
The shifted Dolbeault complex $\OO^{0,\bullet-1}(\pT)$ of a K\"ahler manifold
$X$  is an $\Linf$ algebra.
The $n$-th multibracket
\[ \lambda_n:\OO^{0,j_1}(\pT)\otimes\cdots\otimes\OO^{0,j_n}(\pT)\to
\OO^{0,j_1+\cdots+j_n+1}(\pT) \] is the composition of the wedge product
\[ \OO^{0,j_1}(\pT)\otimes\cdots\otimes\OO^{0,j_n}(\pT)\to
\OO^{0,j_1+\cdots+j_n}(\otimes^n \pT) \] with the map
\[ \OO^{0,j_1+\cdots+j_n}(\otimes^n \pT) \to
\OO^{0,j_1+\cdots+j_n+1}(\pT) \]
associated to $R_n\in\OO^{0,1}\big(\Hom(\otimes^n \pT,\pT)\big)$ in the obvious way.
\end{corollary}

\subsection{Proofs}\label{Subsectin:proofsofmaintheorems}

This part is devoted to the proofs of the theorems claimed  earlier
in  this section.
For convenience, we set $\widetilde{b}=j(b)$,
$\forall b\in\Gamma(B)$.

\subsubsection{Atiyah class as a Lie bracket}

Below we  follow the notations as introduced in Section \ref{section:341}.

\begin{lemma}\label{Lemma6}
For any $a_1,a_2\in\sections{A}$ and $b\in\sections{B}$, we have
\[ \lie{\anadelta_b a_1}{a_2}+\lie{a_1}{\anadelta_b a_2}
-\anadelta_b\lie{a_1}{a_2}=\anadelta_{\nabla_{a_1}b}a_2-\anadelta_{\nabla_{a_2}b}a_1 .\]
\end{lemma}

\begin{proof}
We have
\begin{align*}
& p\big( \lie{\lie{\tb}{a_1}}{a_2}+\lie{\lie{a_1}{a_2}}{\tb}+\lie{\lie{a_2}{\tb}}{a_1} \big) \\
=& p\lie{p\lie{\tb}{a_1}+\widetilde{q\lie{\tb}{a_1}}}{a_2}-p\lie{\tb}{\lie{a_1}{a_2}}
+p\lie{p\lie{a_2}{\tb}+\widetilde{q\lie{a_2}{\tb}}}{a_1} \\
=& \lie{p\lie{\tb}{a_1}}{a_2}-p\lie{\widetilde{q\lie{a_1}{\tb}}}{a_2}-p\lie{\tb}{\lie{a_1}{a_2}}
+\lie{a_1}{p\lie{\tb}{a_2}}+p\lie{\widetilde{q\lie{a_2}{\tb}}}{a_1} \\
=& \lie{\anadelta_b a_1}{a_2}-\anadelta_{\nabla_{a_1}b}a_2-\anadelta_b\lie{a_1}{a_2}
+\lie{a_1}{\anadelta_b a_2}+\anadelta_{\nabla_{a_2}b}a_1
.\end{align*}
The result follows from the Jacobi identity of the Lie algebroid $L$.
\end{proof}

Note that a bundle map $\omega:(\wedge^k A)\otimes(\otimes^l B)\to E$
determines a bundle map
\begin{equation}
\label{eq:ms}
\cev{\omega}:\wedge^k A\to (\otimes^l B^*)\otimes E
\end{equation}
and vice versa.

\begin{lemma}\label{Lemma9}
For any bundle map $\omega:(\wedge^k A)\otimes(\otimes^l B)\to E$, any $a_0,\dots,a_k\in\sections{A}$
and any $b_0,\dots,b_l\in\sections{B}$, we have
\begin{align*}
& (-1)^{k} \big(\partial^A\parnab\omega+\parnab\partial^A\omega\big)(a_0,\cdots,a_k;b_0,\dots,b_l)  \\
& \quad = \sum_{i=0}^k (-1)^{i}  \left\langle \nabla_{a_i}\nabla_{\widetilde{b_0}}\big(\cev{\omega}(\widehat{a_i})\big)
-\nabla_{\widetilde{b_0}}\nabla_{a_i}\big(\cev{\omega}(\widehat{a_i})\big)
-\nabla_{\lie{a_i}{\widetilde{b_0}}}\big(\cev{\omega}(\widehat{a_i})\big) \middle| b_1\otimes\cdots\otimes b_l \right\rangle  \\
& \quad = \sum_{i=0}^k (-1)^{i} \left\{ R^\nabla_E(a_i;b_0)\cdot\omega(\widehat{a_i};\widehat{b_0})
 - \sum_{j=1}^l \omega\big(\widehat{a_i};b_1,\cdots,R^\nabla_B(a_i;b_0)\cdot b_j,\cdots,b_l\big) \right\}
,\end{align*}
where $\widehat{a_i}$ stands for $a_0\wedge\cdots\wedge a_{i-1}\wedge a_{i+1}\wedge\cdots\wedge a_k$
and $\widehat{b_0}$ for $b_1\otimes\cdots\otimes b_l$.
\end{lemma}

\begin{proof}[Sketch of proof]
The first equality follows from a cumbersome computation
at the last step of which use is made of Lemma~\ref{Lemma6}.
The second equality is immediate.
\end{proof}

Given $\mu\in\sections{(\wedge^{k_1}A^*)\otimes(\otimes^{l_1}B^*)\otimes B}$,
$\nu\in\sections{(\wedge^{k_2}A^*)\otimes(\otimes^{l_2}B^*)\otimes B}$,
and arbitrary sections $b_1,\dots,b_{l_1},b'_1,\dots,b'_{l_2}$ of $B$,
let \[
 \creep{\mu\big(b_1,\cdots,b_{i-1},\nu(b'_1,\cdots,b'_{l_2}),b_{i+1},
\cdots,b_{l_1}\big)} \in \sections{ \wedge^{k_1+k_2}A^* \otimes B} ,\]
be defined, as   a bundle map $\wedge^{k_1+k_2}A\to B$,
sending $a_1\wedge\cdots\wedge a_{k_1+k_2}$  to
\begin{multline*} \sum_{\sigma\in\shuffle{k_1}{k_2}} \sgn(\sigma)
\mu\big(a_{\sigma(1)},\cdots,a_{\sigma(k_1)};b_1,\cdots,b_{i-1}, \\
\nu( a_{\sigma(k_1+1)},\cdots,a_{\sigma(k_1+k_2)};b'_1,\cdots,b'_{l_2})
,b_{i+1},\cdots,b_{l_1}\big) .\end{multline*}
In particular, if $\mu=\alpha_1\otimes\beta_1\otimes u$ and $\nu=\alpha_2\otimes\beta_2\otimes v$
with $\alpha_1\in\sections{\wedge^{k_1}A^*}$, $\alpha_2\in\sections{\wedge^{k_2}A^*}$,
$\beta_1\in\sections{\otimes^{l_1}B^*}$, $\beta_2\in\sections{\otimes^{l_2}B^*}$,
and $u,v\in\sections{B}$, then
\begin{multline*}
\creep{\mu\big(b_1,\cdots,b_{i-1},\nu(b'_1,\cdots,b'_{l_2}),b_{i+1},\cdots,b_{l_1}\big)} = \\
\beta_1(b_1,\cdots,b_{i-1},v,b_{i+1},\cdots,b_{l_1})\beta_2(b'_1,\cdots,b'_{l_2}) \cdot(\alpha_1\wedge\alpha_2)\otimes u
.\end{multline*}

\begin{corollary}\label{Corollary10}
For any $n\geq 2$ and $b_0,\dots,b_n\in\sections{B}$, we have
\begin{multline*}
- \big((\partial^A\parnab+\parnab\partial^A)R_n\big)(b_0,\cdots,b_n) = \creep{R_2\big(b_0,R_n(b_1,\cdots,b_n)\big)} \\
+\sum_{j=1}^n \creep{R_n\big(b_1,\cdots, b_{j-1}, R_2(b_0,b_j), b_{j+1},
\cdots,b_n\big)}
.\end{multline*}
\end{corollary}

\begin{proof}
Apply Lemma~\ref{Lemma9} to $\omega=R_n$.
\end{proof}

\begin{corollary}\label{Corollary11}
For any $b_0,b_1,b_2\in\sections{B}$, we have
\begin{multline}
- \big(\partial^A R_3\big)(b_0,b_1,b_2) = \creep{R_2\big(b_0,R_2(b_1,b_2)\big)}
+\creep{R_2\big(R_2(b_0,b_1),b_2\big)} \\
 +\creep{R_2\big(b_1,R_2(b_0,b_2)\big)}
.
\label{eq:tiger}
\end{multline}
\end{corollary}

\begin{proof}
Since $\partial^A R_2=0$ and $\parnab R_2=R_3$, taking $n=2$ in Corollary~\ref{Corollary10} yields the result.
\end{proof}

\begin{proof}[Sketch of proof of Theorem~\ref{laoidc}]
The interchange isomorphism $\tau:B[-1]\otimes B[-1]\to B[-1]\otimes B[-1]$
 is the image in $D^+ (\category{A})$
of the chain map $\tau:B[-1]\otimes B[-1]\to B[-1]\otimes B[-1]$ given by $\tau(b_1\otimes b_2)=- b_2\otimes b_1$,
$\forall b_1,b_2\in B$ --- the negative sign is due to
 $B[-1]$ being a complex concentrated in degree 1 (see
Equation~\eqref{braiding}).
Recall that $R_2$ is a cocycle (w.r.t.\ $\partial^A$).
Its cohomology class $\alpha_B$, the Atiyah class of $B$, can be seen as an element of
$\Hom_{D^+(\category{A})}(B[-1]\otimes B[-1],B[-1])$.
Proposition~\ref{Lemma1} implies the equality $\alpha_B\rond\tau=-\alpha_B$
in $\Hom_{D^+(\category{A})}(B[-1]\otimes B[-1],B[-1])$.
Corollary~\ref{Corollary11} implies that the Jacobi identity
$\alpha_B\rond(\id\varotimes\alpha_B)=\alpha_B\rond(\alpha_B\varotimes\id)
+\alpha_B\rond(\id\varotimes\alpha_B)\rond(\tau\varotimes\id)$
holds in $D^+(\category{A})$. Indeed, each of the terms
in the right hand side of Equation \eqref{eq:tiger}
can be interpreted as a Yoneda product, a composition of morphisms in the derived category.
\end{proof}

\subsubsection{Jacobi identity up to homotopy}

Consider the cochain complex $(V[-1],\partial^A)$,
where the graded vector space $V=\bigoplus_{k=0}^{\infty} V_k$
is given by $V_k=\sections{\wedge^k A^*\otimes B}$ so that,
if $\xi\in\sections{\wedge^k A^*}$ and $b\in\sections{B}$,
then $\xi\otimes b\in \big(V[-1]\big)_{k+1}$.

\begin{lemma}\label{Lemma12}
The graded linear map $\lambda:V[-1]\varotimes V[-1]\to V[-1]$ given by
\[ \lambda(v_1\varotimes v_2) = (-1)^{k_2} \xi_1\wedge\xi_2\wedge R_2(b_1,b_2) \]
for any $v_1=\xi_1\otimes b_1\in (V[-1])_{k_1+1}$ and $v_2=\xi_2\otimes b_2\in (V[-1])_{k_2+1}$ is a chain map.
\end{lemma}

\begin{proof}
A straightforward computation yields
\[ \big(\partial^A\rond\lambda-\lambda\rond\partial^A\big)(v_1\varotimes v_2)
=(-1)^{k_1}\xi_1\wedge\xi_2\wedge\big(\partial^A R_2\big)(b_1,b_2) .\]
The result follows from $\partial^A R_2=0$ (see Theorem~\ref{Thm:atiyahclass} and the definition \eqref{r2} of $R_2$).
\end{proof}

Now, consider the interchange isomorphism
$\tau:V[-1]\varotimes V[-1]\to V[-1]\varotimes V[-1]$ given by
$\tau(v_1\varotimes v_2)=(-1)^{\abs{v_1}\abs{v_2}}v_2\varotimes v_1$.

\begin{lemma}\label{Lemma13}
The chain map $\lambda$ is skew-symmetric up to  homotopy:
 \[ \lambda+\lambda\rond\tau=\partial^A\rond\Theta+\Theta\rond\partial^A ,\]
where the graded map $\Theta:V[-1]\varotimes V[-1]\to V[-2]$ is given by
\[ \Theta(v_1\varotimes v_2) = (-1)^{k_1} \xi_1\wedge\xi_2\otimes\beta(b_1,b_2) \]
for any $v_1=\xi_1\otimes b_1\in (V[-1])_{k_1+1}$ and $v_2=\xi_2\otimes b_2\in (V[-1])_{k_2+1}$.
\end{lemma}

\begin{proof}
Straightforward computations yield
\[ \big(\lambda+\lambda\rond\tau\big)(v_1\varotimes v_2)
=(-1)^{k_2}\xi_1\wedge\xi_2\wedge\big(R_2(b_1,b_2)-R_2(b_2,b_1)\big) \]
and
\[ \big(\partial^A\rond\Theta+\Theta\rond\partial^A\big)(v_1\varotimes v_2)
=(-1)^{k_2}\xi_1\wedge\xi_2\wedge  \big( \big(\partial^A\beta\big)(b_1,b_2)
\big)  \]
The result follows from Proposition~\ref{Lemma1}.
\end{proof}

\begin{lemma}\label{Lemma14}
The chain map $\lambda$ satisfies the Jacobi identity up to a homotopy:
\[ -\lambda\rond(\id\varotimes\lambda)+\lambda\rond(\lambda\varotimes\id)
+\lambda\rond(\id\varotimes\lambda)\rond(\tau\varotimes\id)
= \partial^A\rond\Xi+\Xi\rond\partial^A ,\]
where the graded map
$\Xi:V[-1]\varotimes V[-1]\varotimes V[-1]\to V[-2]$ is given by
\[ \Xi(v_0\varotimes v_1\varotimes v_2) = (-1)^{k_0+k_2} \xi_0\wedge\xi_1\wedge\xi_2\wedge R_3(b_0,b_1,b_2) \]
for any $v_i=\xi_i\otimes b_i\in (V[-1])_{k_i+1}$ with $i\in\{0,1,2\}$.
\newline
\end{lemma}

\begin{proof}
Straightforward computations yield
\begin{gather*}
\big( \lambda\rond(\id\varotimes\lambda) \big)(v_0\varotimes v_1\varotimes v_2)
= (-1)^{k_1}\xi_0\wedge\xi_1\wedge\xi_2\wedge\creep{R_2(b_0,R_2(b_1,b_2))} ,\\
\big( \lambda\rond(\lambda\varotimes\id) \big)(v_0\varotimes v_1\varotimes v_2)
= -(-1)^{k_1}\xi_0\wedge\xi_1\wedge\xi_2\wedge\creep{R_2(R_2(b_0,b_1),b_2)} ,\\
\big( \lambda\rond(\id\varotimes\lambda)\rond(\tau\varotimes\id) \big)(v_0\varotimes v_1\varotimes v_2)
= -(-1)^{k_1}\xi_0\wedge\xi_1\wedge\xi_2\wedge\creep{R_2(b_1,R_2(b_0,b_2))} ,\\
\intertext{and}
\big( \partial^A\rond\Theta+\Theta\rond\partial^A \big)(v_0\varotimes v_1\varotimes v_2)
= (-1)^{k_1}\xi_0\wedge\xi_1\wedge\xi_2\wedge \big( \big(\partial^A R_3\big)(b_0,b_1,b_2) \big)
.\end{gather*}
The result follows from Corollary~\ref{Corollary11}.
\end{proof}

Theorem~\ref{Thm:CsuperLiealgebraobject} immediately follows from
Lemmas~\ref{Lemma12}, \ref{Lemma13}, and~\ref{Lemma14}.

Note that Theorem~\ref{Thm:CsuperLiealgebraobject} could also be seen
as a corollary of Theorems~\ref{Thm:CisSHLeiniz} and~\ref{thm:module}.

\subsubsection{\leibnizinfinityone\ algebra (and modules) arising from a Lie pair}

\begin{lemma}\label{Lemma15}
For any $n\geq 3$ and any arbitrary sections $b_1,\dots,b_n$ of $B$, we have
\begin{multline*}
-\big(\partial^A R_n\big)(b_1,\dots,b_n)=
\sum_{\substack{i+j=n+1 \\ i\geq 2 \\ j\geq 2}}\sum_{k=j}^n\sum_{\sigma\in\shuffle{k-j}{j-1}} \\
\creep{R_i\big(b_{\sigma(1)},\cdots,b_{\sigma(k-j)},R_j(b_{\sigma(k+1-j)},\cdots,
b_{\sigma(k-1)},b_k),b_{k+1},\cdots,b_n\big)}
.\end{multline*}
\end{lemma}

\begin{proof}
We reason by induction.  The formula holds for $n=3$ by Corollary~\ref{Corollary11}.
Assuming the formula holds for $n=N$, we get
\begin{multline*}
\big(\parnab\partial^A R_N\big)(b_0,\dots,b_N)=\big(\nabla_{b_0}(\partial^A R_N)\big)(b_1,\dots,b_N)
= \sum_{\substack{i+j=N \\ i\geq 2 \\ j\geq 2}}\sum_{k=j}^N\sum_{\sigma\in\shuffle{k-j}{j-1}} \\
\Big\{ \creep{R_{i+1}\big(b_0,b_{\sigma(1)},\cdots,b_{\sigma(k-j)},R_j(b_{\sigma(k+1-j)},
\cdots,b_{\sigma(k-1)},b_k),b_{k+1},\cdots,b_n\big)} \\
+ \creep{R_i\big(b_{\sigma(1)},\cdots,b_{\sigma(k-j)},R_{j+1}(b_0,b_{\sigma(k+1-j)},\cdots,
b_{\sigma(k-1)},b_k),b_{k+1},\cdots,b_n\big)} \Big\}
.\end{multline*}
Observing that
\[ \partial^A R_{N+1}=(\partial^A\parnab+\parnab\partial^A) R_N - \parnab\partial^A R_N \]
and recalling Corollary~\ref{Corollary10},
it is easy to check that the desired formula holds for $n=N+1$ as well.
\end{proof}

\begin{lemma}\label{Lemma16}
For any bundle map $\omega:(\wedge^k A)\otimes(\otimes^l B)\to B$ and any $b_1,\dots,b_l\in\sections{B}$, we have
\[  \partial^A\big(\cev{\omega}(b_1,\dots,b_l)\big)
-\big(\cev{\partial^A\omega}\big)(b_1,\dots,b_l)=
(-1)^k\sum_{j=0}^l \creep{\omega(b_1,\cdots, b_{j-1},
\partial^A b_j, b_{j+1}, \cdots,b_l)} ,\]
where $\cev{\omega}$ is defined by Equation \eqref{eq:ms}.
\end{lemma}

\begin{proof}
For any $a_0,\dots,a_k\in\sections{A}$, we have
\begin{align*}
&\left\langle \partial^A\big(\cev{\omega}  (b_1,\dots,b_l)\big)-
\big(\cev{\partial^A\omega}\big)(b_1,\dots,b_l)
 \middle| a_0\wedge\cdots\wedge a_k \right\rangle \\
=\; & \sum_{j=0}^l \sum_{i=0}^k (-1)^i \omega( a_0,\cdots,a_{i-1},a_{i+1},\cdots,
a_k ;b_1,\cdots,\nabla_{a_i}b_j,\cdots,b_l) \\
=\; & (-1)^k \sum_{j=0}^l \sum_{\sigma\in\shuffle{k}{1}}\sgn(\sigma)\omega( a_{\sigma(0)},\cdots,
a_{\sigma(k-1)}; b_1,\cdots,\nabla_{a_{\sigma(k)}}b_j,\cdots,b_l) \\
=\; & (-1)^k \sum_{j=0}^l \left\langle \creep{\omega(b_1,\cdots,
b_{j-1}, \partial^A b_j, b_{j+1},  \cdots,b_l)} \middle|
a_0\wedge\cdots\wedge a_k \right\rangle . \qedhere \end{align*}
\end{proof}

\begin{proof}[Proof of Theorem~\ref{Thm:CisSHLeiniz}]
We only need to check that the generalized Leibniz
 identity \eqref{Eqt:generalJacobi} holds.
Since $\lambda_1=\partial^A$ and $(\partial^A)^2=0$, Equation~\eqref{Eqt:generalJacobi} is obviously true for $n=1$.
Let $n\geq 2$ and $v_i=\xi_i\otimes b_i\in\sections{\wedge^{p_i}A^*\otimes B}$ for all $i\in\{1,\dots,n\}$.
The l.h.s.\ of  \eqref{Eqt:generalJacobi} is
\begin{multline*}
\sum_{1\le j\le k\le n}\sum_{\sigma\in\shuffle{k-j}{j-1}}
\Koszul{\sigma;v_1,\cdots,v_{k-1}}
(-1)^{\abs{v_{\sigma(1)}}+\abs{v_{\sigma(2)}}+\cdots+\abs{v_{\sigma(k-j)}}} \\
\lambda_{n-j+1}\big(v_{\sigma(1)},\cdots,v_{\sigma(k-j)},
\lambda_j(v_{\sigma(k+1-j)},\cdots,v_{\sigma(k-1)},v_k),
v_{k+1},\cdots, v_{n}\big) .
\end{multline*}
Separating the terms involving $\lambda_1$ (aka $\partial^A$) from the others, it can be rewritten as
\begin{multline*}
\partial^A\big(\lambda_n(v_1,\cdots,v_n)\big)
+ \sum_{\substack{i+j=n+1 \\ i\geq 2 \\ j\geq 2}}\sum_{k=j}^n\sum_{\sigma\in\shuffle{k-j}{j-1}}
\Koszul{\sigma;\xi_1,\cdots,\xi_{k-1}}(-1)^{p_{\sigma(1)}+\cdots+p_{\sigma(k-j)}} \\
\lambda_{i}\big(v_{\sigma(1)},\cdots,v_{\sigma(k-j)},\lambda_j(v_{\sigma(k+1-j)},\cdots,
v_{\sigma(k-1)},v_k),v_{k+1},\cdots, v_{n}\big) \\
+\sum_{k=1}^n (-1)^{p_1+\cdots+p_{k-1}} \lambda_n\big(v_1,\cdots,v_{k-1},
(\partial^A \xi_k)\otimes b_k+(-1)^{p_k}\xi_k\otimes(\partial^A b_k) ,v_{k+1},\cdots,v_n\big).
\end{multline*}
Since  each $\lambda_k$  is given by
 Equation ~\eqref{eqn:temp1} in terms of the corresponding $R_k$, it in turn
becomes
\begin{multline*}
\partial^A\big((-1)^{p_1+\cdots+p_n}\xi_1\wedge\cdots\wedge\xi_n\wedge R_n(b_1,\cdots,b_n)\big) \\
+ \sum_{\substack{i+j=n+1 \\ i\geq 2 \\ j\geq 2}}\sum_{k=j}^n\sum_{\sigma\in\shuffle{k-j}{j-1}}
\Koszul{\sigma;\xi_1,\cdots,\xi_{k-1}}
\xi_{\sigma(1)}\wedge\cdots\wedge\xi_{\sigma(k-1)}\wedge\xi_k\wedge\cdots\wedge\xi_n\wedge \\
\creep{R_i\big(b_{\sigma(1)},\cdots,b_{\sigma(k-j)},R_j(b_{\sigma(k+1-j)},\cdots,
b_{\sigma(k-1)},b_k),b_{k+1},\cdots,b_{n}\big) } \\
+ \sum_{k=1}^n (-1)^{1+p_k+p_{k+1}+\cdots+p_n}\xi_1\wedge\cdots\wedge\xi_{k-1}\wedge
\partial^A\xi_k\wedge\xi_{k+1}\wedge\cdots\wedge\xi_n\wedge R_n(b_1,\cdots,b_n) \\
+ \sum_{k=1}^n \xi_1\wedge\cdots\wedge\xi_n\wedge
\creep{R_n(b_1,\cdots, b_{k-1}, \partial^A b_k, b_{k+1}, \cdots,b_n)}
,\end{multline*}
which simplifies to
\begin{multline*}
\xi_1\wedge\cdots\wedge\xi_n\wedge \Big\{ \partial^A\big(R_n(b_1,\cdots,b_n)\big)
+\sum_{\substack{i+j=n+1 \\ i\geq 2 \\ j\geq 2}}\sum_{k=j}^n\sum_{\sigma\in\shuffle{k-j}{j-1}} \\
\creep{R_i\big(b_{\sigma(1)},\cdots,b_{\sigma(k-j)},R_j(b_{\sigma(k+1-j)},\cdots,
b_{\sigma(k-1)},b_k),b_{k+1},\cdots,b_{n}\big) } \\
+ \sum_{k=1}^n \creep{R_n(b_1,\cdots, b_{k-1}, \partial^A b_k, b_{k+1}, \cdots,b_n)} \Big\}
.\end{multline*}
The result now follows from Lemmas~\ref{Lemma15} and~\ref{Lemma16}.
\end{proof}

The proofs of Theorem~\ref{thm:module}, Theorem~\ref{thm:mainC}
and Corollary~\ref{cor:Lie} go along the same line \emph{mutatis mutandis}.

\subsubsection{$L_\infty$ rather than \leibnizinfinity}

\begin{lemma}\label{Lemma2}
For any $b_0,b_1,b_2\in\sections{B}$, we have
\[ R_3(b_0,b_1,b_2)-R_3(b_1,b_0,b_2)=R_2\big(\beta(b_0,b_1),b_2\big)-
\big(\partial^A\Omega\big)(b_0,b_1)\cdot b_2 ,\]
where $\beta:\wedge^2 B\to B$ and
$\Omega:\wedge^2 B\to \End B$ are bundle maps as  introduced
in Equations \eqref{eq:snake1}-\eqref{eq:snake2}.
\end{lemma}

\begin{proof}
The second  Bianchi identity $d^\nabla R^\nabla=0$ holds for the curvature
$R^\nabla:\wedge^2 L\to\End B$ of the $L$-connection $\nabla$ on $B$
extending the $A$-action.
Hence, for any $a\in\sections{A}$ and $b_0,b_1,b_2\in\sections{B}$, we have
\begin{align*}
0=& \big(d^\nabla R^\nabla\big)(a,\tb_0,\tb_1) \\
=& \nabla_a\big( R^\nabla(\tb_0,\tb_1) \big) - \nabla_{\tb_0}\big( R^\nabla(a,\tb_1) \big)
+ \nabla_{\tb_1}\big( R^\nabla(a,\tb_0) \big) \\
& - R^\nabla(\lie{a}{\tb_0},\tb_1) + R^\nabla(\lie{a}{\tb_1},\tb_0)
- R^\nabla(\lie{\tb_0}{\tb_1},a) \\
=& \nabla_a\big( R^\nabla(\tb_0,\tb_1) \big) - \nabla_{\tb_0}\big( R^\nabla(a,\tb_1) \big)
+ \nabla_{\tb_1}\big( R^\nabla(a,\tb_0) \big) \\
& - R^\nabla(\widetilde{\nabla_a b_0},\tb_1) + R^\nabla(\anadelta_{b_0}a,\tb_1)
+ R^\nabla(\widetilde{\nabla_a b_1},\tb_0) - R^\nabla(\anadelta_{b_1}a,\tb_0) \\
& - R^\nabla\big(\alpha(b_0,b_1),a\big) - R^\nabla(\widetilde{\nabla_{\tb_0}b_1},a)
+ R^\nabla(\widetilde{\nabla_{\tb_1}b_0},a) +R^\nabla\big(\widetilde{\beta(b_0,b_1)},a\big)
\end{align*}
and thus
\begin{align*}
0=& \big(d^\nabla R^\nabla\big)(a,\tb_0,\tb_1) \cdot b_2 \\
=& \big(\parnab R^\nabla_B\big)(a;b_0,b_1) \cdot b_2
- \big(\parnab R^\nabla_B\big)(a;b_1,b_0) \cdot b_2
- R^\nabla_B\big(a,\beta(b_0,b_1)\big) \cdot b_2 \\
& + \nabla_a\big(R^\nabla(\tb_0,\tb_1)\big)\cdot b_2
- R^\nabla\big(\widetilde{\nabla_a b_0},\tb_1\big) \cdot b_2
- R^\nabla\big(\tb_0,\widetilde{\nabla_a b_1}\big) \cdot b_2
\end{align*}
or, equivalently,
\[ 0=R_3(b_0,b_1,b_2)-R_3(b_1,b_0,b_2)-R_2\big(\beta(b_0,b_1),b_2\big)
+\big(\partial^A\Omega\big)(b_0,b_1)\cdot b_2 .\qedhere\]
\end{proof}

\begin{lemma}\label{Lemma3}
For any $a\in\sections{A}$ and $b_0,b_1\in\sections{B}$, we have
\[ \lie{\alpha(b_0,b_1)}{a}+\alpha(\nabla_a b_0,b_1)+\alpha(b_0,\nabla_a b_1)
= \anadelta_{b_0}\anadelta_{b_1}a-\anadelta_{b_1}\anadelta_{b_0}a
-\anadelta_{\widetilde{q\lie{\tb_0}{\tb_1}}} a .\]
\end{lemma}

\begin{proof}
We have
\begin{align*}
& p\big( \lie{\tb_1}{\lie{\tb_0}{a}}+\lie{\tb_0}{\lie{a}{\tb_1}}+\lie{a}{\lie{\tb_1}{\tb_0}} \big) \\
=& p\lie{\tb_1}{p\lie{\tb_0}{a}}+p\lie{\tb_0}{p\lie{a}{\tb_1}}+p\lie{a}{p\lie{\tb_1}{\tb_0}} \\
& +p\lie{\tb_1}{\widetilde{q\lie{\tb_0}{a}}}+p\lie{\tb_0}{\widetilde{q\lie{a}{\tb_1}}}
+p\lie{a}{\widetilde{q\lie{\tb_1}{\tb_0}}} \\
=& \anadelta_{b_1}\anadelta_{b_0}a-\anadelta_{b_0}\anadelta_{b_1}a
+ p\lie{a}{\alpha(b_1,b_0)}+p\lie{\widetilde{\nabla_a b_0}}{\tb_1}
+p\lie{\tb_0}{\widetilde{\nabla_a b_1}}+p\lie{\widetilde{q\lie{\tb_0}{\tb_1}}}{a} \\
=&  \anadelta_{b_1}\anadelta_{b_0}a-\anadelta_{b_0}\anadelta_{b_1}a
+ \lie{\alpha(b_0,b_1)}{a} +\alpha(\nabla_a b_0,b_1)+\alpha(b_0,\nabla_a b_1)
+\anadelta_{\widetilde{q\lie{\tb_0}{\tb_1}}}a
.\end{align*}
The result follows from the Jacobi identity of $L$.
\end{proof}

\begin{lemma}\label{Lemma4}
For any $n\geq 3$, $a\in\sections{A}$ and $b_0,b_1,\dots,b_n\in\sections{B}$, we have
\begin{multline*}
R_{n+1}(a;b_0,b_1,b_2,\cdots,b_n)-R_{n+1}(a;b_1,b_0,b_2,\cdots,b_n) = \\
\Omega(b_0,b_1)\cdot R_{n-1}(a;b_2,\cdots,b_n)
-\sum_{j=2}^n R_{n-1}\big(a;b_2,\cdots,\Omega(b_0,b_1)\cdot b_j,\cdots,b_n\big) \\
+\nabla_{\alpha(b_0,b_1)}\big(R_{n-1}(a;b_2,\cdots,b_n)\big)
-\sum_{j=2}^n R_{n-1}\big(a;b_2,\cdots,\nabla_{\alpha(b_0,b_1)} b_j,\cdots,b_n\big) \\
-R_{n-1}\big( \lie{\alpha(b_0,b_1)}{a}+\alpha(\nabla_a b_0,b_1)+\alpha(b_0,\nabla_a b_1);b_2,\cdots,b_n \big) \\
+R_n\big(a;\beta(b_0,b_1),b_2,\cdots,b_n\big)
.\end{multline*}
Here $R_{k}(a;b_1,\cdots,b_k)$ means $R_{k}(b_1,\cdots,b_k)(a)$.
\end{lemma}

\begin{proof}[Sketch of proof]
Straightforward computation at the last step of which use is made of Lemma~\ref{Lemma3}.
\end{proof}

\begin{proposition}\label{Proposition5}
Let $L=A\bowtie B$ be a matched pair of Lie algebroids endowed with a flat torsion free $B$-connection on $B$.
(These data determine a splitting of the short exact sequence of vector bundles
$0\to A\to L\to B\to 0$
and an $L$-connection on $B$ extending the $A$-action
such that the three associated bundle maps $\alpha$, $\beta$, and $\Omega$ are trivial.)
Then each $R_n:\otimes^n B\to\Hom(A,B)$ is totally symmetric in its $n$ arguments.
\end{proposition}

\begin{proof}
It follows from Proposition~\ref{Lemma1} and Lemma~\ref{Lemma2}
that $R_2$ and $R_3$ are invariant under the permutation of their first two arguments.
By Lemma~\ref{Lemma4}, the same property holds for all higher $R_n$.
Moreover, it is easy to see that, if $R_n$ is symmetric in its $n$ arguments,
then $R_{n+1}$ is symmetric in its last $n$ arguments since $R_{n+1}=\parnab R_n$.
The result follows by induction.
\end{proof}

Theorem~\ref{Thm:conditionforCtobeSHLie}, which says that $V[-1]$ is an
$L_\infty$-algebra when the assumptions of Proposition~\ref{Proposition5} are satisfied,
is an immediate consequence of Proposition~\ref{Proposition5}
and Theorem~\ref{Thm:CisSHLeiniz}.


\bibliographystyle{amsplain}
\bibliography{biblio}
\end{document}